# A Theory of Structural Independence

**Matthias G. Mayer**
matthias.georg.mayer@gmail.com

15th April 2025

## Abstract

*Structural independence* is the (conditional) independence that arises from the structure rather than the precise numerical values of a distribution. We develop this concept and relate it to $d$-separation and structural causal models.

Formally, let $U = (U_i)_{i \in I}$ be an independent family of random elements on a probability space $(\Omega, \mathcal{A}, \mathbb{P})$. Let $X, Y$ and $Z$ be arbitrary $\sigma(U)$-measurable random elements. We characterize all independences $X \perp\!\!\!\perp Y \mid Z$ implied by the independence of $U$ and call these independences *structural*. Formally, these are the independences which hold in all probability measures $P$ that render $U$ independent and are absolutely continuous with respect to $\mathbb{P}$, i.e. for all such $P$, it needs to hold that $X \perp\!\!\!\perp_P Y \mid Z$.

We introduce the history $\mathcal{H}(X|Z) : \Omega \to \mathfrak{P}(I)$, a combinatorial object that measures the dependence of $X$ on $U_i$ for each $i \in I$ given $Z$. The independence of $X$ and $Y$ given $Z$ is implied by the independence of $U$ if and only if $\mathcal{H}(X|Z) \cap \mathcal{H}(Y|Z) \overset{\text{a.s.}}{=} \emptyset$ with respect to $\mathbb{P}$.

Finally, we apply this $d$-separation like criterion in structural causal models to discover a causal direction in a toy setting.

# 1 Introduction

The purpose of this work[1] is to investigate when the independence of random elements is implied by the independence of a family of random elements. We give a review of basic, well-known definitions and theorems in Appendix B that will also be used in the rest of this introduction. In Section 2, we review $d$-separation, what it is used for, and how it relates to the phenomenon we study.

Let $I$ be an arbitrary index set. Let $(\Omega, \mathcal{A}, \mathbb{P})$ be a probability space. Let $(U_i)_{i \in I}$ be an independent family of random elements.

**In the following** $X, Y$ and $Z$ will always refer to $\sigma(U)$-measurable random elements on $(\Omega, \mathcal{A})$.

We investigate the question: "Which independences in $\mathbb{P}$ are implied by $U_i$ being independent?" To formalize this, we let

$$\triangle^{\times} := \Big\{ P : \mathcal{A} \to \mathbb{R} \mid P \text{ is a probability measure, } (U_i)_{i \in I} \text{ is an independent family w.r.t. } P \text{ and } P \sim \mathbb{P} \Big\}$$

and ask when $\forall P \in \triangle^{\times} : X \perp\!\!\!\perp_P Y \mid Z$ holds. More specifically, we want to characterize this statement uniformly over all choices of $X, Y$ and $Z$ without quantifying over $\triangle^{\times}$.

**Remark 1.1**: The requirement that $P \sim \mathbb{P}$ can be relaxed to $P \ll \mathbb{P}$ without changing the resulting characterization. It is however convenient for calculations to require $P \sim \mathbb{P}$, s.t. $\frac{\mathrm{d}P}{\mathrm{d}\mathbb{P}}$ can be chosen to be positive. If an independence holds for all $P \sim \mathbb{P}$, the case $P \ll \mathbb{P}$ follows as a corollary by density (Lemma A.9).

In future work, we might investigate the case where we remove the reference measure and the condition $P \sim \mathbb{P}$ from $\triangle^{\times}$, c.f. Section 10.6.

From standard probability we can immediately conclude a collection of unconditional independences. For $J \subseteq I$, let $U_J$ denote the family $(U_i)_{i \in J}$.

**Lemma 1.2**: Let $J, K \subseteq I$ be disjoint. Let $\sigma(X) \subseteq \sigma(U_J)$ and $\sigma(Y) \subseteq \sigma(U_K)$. Then $X \perp\!\!\!\perp_P Y$ for all $P \in \triangle^{\times}$.

*Proof*: Let $A \in \sigma(X)$ and $B \in \sigma(Y)$. By assumption $A \in \sigma(U_J)$ and $B \in \sigma(U_K)$. $A \perp\!\!\!\perp_P B$ follows immediately because $U_J \perp\!\!\!\perp_P U_K$, since $U$ is an independent family w.r.t. $P$, and $J \cap K = \emptyset$. □

However, it is a priori not clear that these are the only independences that follow. We want to see, whether $\forall P \in \triangle^{\times} : X \perp\!\!\!\perp_P Y$ implies that there are $J, K \subseteq I$ disjoint, s.t. $\sigma(X) \subseteq \sigma(U_J)$ and $\sigma(Y) \subseteq \sigma(U_K)$. We will see that this is indeed the case in Theorem 6.4.1

[1] A version of this work was submitted as a Master's thesis at the University of Bonn.

The interesting theory begins once we condition on a random element $Z$. When do we have the conditional independence $X \perp\!\!\!\perp_P Y \mid Z$ for all $P \in \triangle^{\times}$? To answer this question, we want to generalize from Lemma 1.2, where we began to characterize this with dependence on disjoint sets of background variables. Conditional on $Z$, we want to carry over this idea and allow this dependence to vary with $Z$. Suppose that $\{Z = z\}$ is not a $P$-nullset. For each $z \in \mathrm{Val}(Z)$, we want to assign $J_z, K_z \subseteq I$, s.t. $\sigma\big(X|_{\{Z=z\}}\big) \subseteq \sigma\big(U_{J_z}|_{\{Z=z\}}\big)$ and $\sigma\big(Y|_{\{Z=z\}}\big) \subseteq \sigma\big(U_{K_z}|_{\{Z=z\}}\big)$. However, even if $J_z$ and $K_z$ are disjoint, we cannot conclude independence.

**Example 1.3**: Let $I = \{1, 2\}$ and $\Omega = \mathbb{F}_2^I$. For $i \in I$, set $U_i = \pi_i$ and let $Z = U_1 + U_2$ (in $\mathbb{F}_2$). Let $C \coloneqq \{Z = 0\}$. We see that for $\omega \in C$, we have $U_1(\omega) = U_2(\omega)$. Trivially, $\sigma(U_i|_C) \overset{\text{a.s.}}{\subseteq} \sigma(U_i|_C|Z)$, so we can set $J_0 \coloneqq \{1\}$ and $K_0 = \{2\}$. Now $\sigma(U_1|_C) \overset{\text{a.s.}}{\subseteq} \sigma\big(U_{J_0}|_C\big)$ and $\sigma(U_2|_C) \overset{\text{a.s.}}{\subseteq} \sigma\big(U_{K_0}|_C\big)$, but $U_1$ and $U_2$ are not independent given $C$ in any product probability distribution on $\Omega$, since $U_1|_C = U_2|_C$.

Therefore, we need to require a niceness condition on the selection of $J_z$. More specifically, we need the condition $U_{J_z}|_{\{Z=z\}} \perp\!\!\!\perp_P U_{\overline{J_z}}|_{\{Z=z\}} \mid \{Z = z\}$. Clearly, using this condition and the same argumentation as in Lemma 1.2, we can then conclude that $X|_{\{Z=z\}} \perp\!\!\!\perp_P Y|_{\{Z=z\}} \mid \{Z = z\}$. Again, we want to see whether this is the only way to conclude independence. And again, the answer is yes (Theorem 6.4.1).

The finite theory reviewed in Section 3 is enough to show this is the finite case.

Once we condition on an arbitrary random element $Z$, we have that $\{Z = z\}$ is a nullset in general. Recall that conditioning on a random element $Z$ is defined by choosing '$\mathbb{E}(X|Z = z)$' almost surely simultaneously for all $z$ through a Radon-Nikodym derivative to circumvent the nullset problem. In the same spirit, we deal with this problem by choosing all $J_z$ almost surely simultaneously with a random index set (Section 4).

In Appendix A, we recall Kakutani's result about mutual absolute continuity of infinite product probability distributions ([1]), and bring it into a form usable for our purposes.

In Section 5, we introduce the history of $X$ given $Z$, a random index set that measures the dependence of $X$ on $U$ given $Z$, and the desiderata it should fulfill. The most important desiderata is the fundamental theorem of structural independence (Desiderata 5.1.2) The fundamental theorem characterizes all independences that are implied by an independent family through $Z$ dependent random index sets through the history.

In Section 6, we examine we introduce the random index set of irrelevance to prove the fundamental theorem of structural independence.

in Section 7, we study properties of the history and structural independence and show that the history is determined uniquely by some desiderata.

In Section 8, we introduce a counterexample that proves that disintegration is not characterized by rectangular atoms in the general case.

In Section 9, we discuss the application of structural independence to structural causal models.

In Section 10, we discuss directions for further research.

# 2 Related work

Structural independence is strongly related to $d$-separation. Just like this work considers probability distributions that render a family $U$ independent, Pearl's theory of Causality introduces a set of probability distributions that satisfy certain independence constraints. Then $d$-separation is a graphical criterion that characterizes which independences are implied by these constraints. We make this more precise in the following definitions.

In the following let $G = (V, E)$ be a directed acyclic graph, where the nodes $V$ represent random variables. For nodes $X, Y \in V$, we write $X \to Y$ for $(X, Y) \in E$. We write $\mathrm{PA}(X) = \{Y : Y \to X\}$ for the parents of a node $X$ in $G$.

**Definition 2.1**: A probability distribution $P$ is compatible with $G$, if it fulfills the Markov condition. This is the case, if any node $X$ is independent (w.r.t. $P$) of all its non-descendants given $\mathrm{PA}(X)$. Let $\triangle(G)$ denote the set of all probability distributions compatible with $G$.

$G$ is often called a Bayesian network or bayesnet for short. We can now ask the question what independences of nodal variables are implied by the Markov condition. More formally, let $X, Y, Z$ be collections of nodes, interpreted as random variables. When is it the case that $\forall P \in \triangle(G) : X \perp\!\!\!\perp_P Y \mid Z$ ? Pearl gives a nice graphical characterization of this statement in [2], Section 1.2.3, see also [3].

**Definition 2.2**: A walk in a graph is a path in the corresponding undirected graph. More precisely, a walk $w$ (from $w_1$ to $w_n$) is a family of nodes $(w_i)_{i=1}^n$, s.t. for all $i \in \{1, ..., n-1\}$, we have $w_i \to w_{i+1}$ or $w_i \leftarrow w_{i+1}$.

**Definition 2.3**: Let $w = (w_i)_{i=1}^n$ be a walk in $G$. For $1 < i < n$, $w_i$ is a collider in $w$, if $w_{i-1} \to w_i \leftarrow w_{i+1}$.

**Definition 2.4** ($d$-connection, $d$-separation): $X$ and $Y$ are $d$-connected given $Z$ (in $G$), if there is a walk $w = (w_i)_{i=1}^n$ from a node in $X$ to a node in $Y$, s.t. $w_i$ is a collider in $w$ if and only if $w_i \in Z$. $X$ and $Y$ are $d$-separated if they are not $d$-connected. In this case, we write $X \perp_d Y \mid Z$. ($d$ is for directed (graph)).

This graphical criterion, $d$-separation, fully characterizes the independence structure of nodal variables (that is implied by the Markov condition).

**Theorem 2.5** (soundness and completeness of $d$-separation):

$$X \perp_d Y \mid Z \Leftrightarrow \forall P \in \triangle(G) : X \perp\!\!\!\perp_P Y \mid Z.$$

*Proof*: See [2], Theorem 1.2.4. or [3] for the original proof. □

This characterization makes it very useful for the problem of causal discovery, see [4] or [2], Chapter 2, for an introduction. Causal discovery is the problem of inferring the simplest graphs that are compatible with a given probability distribution. By Occam's razor, these then constitute the best guess of what the (graphical) causal structure of the data is. It is called 'causal', because the graphs are imbued with a causal meaning, where the arrows correspond to proximal causality. Informally, a node $X$ is a probabilistic (noisy) function of its parents, where the noise is interpreted as all the factors in the world that are not included in this particular model. This can be made formal by structural (or functional) causal models, see [2], Section 1.4. and Section 9.

Classical causal discovery, as in [2], Chapter 2, can simplified be described as follows. Given a distribution $P$ on a measurable space with random variables $V$, extract all (conditional) independences between nodes and interpret them a the $d$-separation relation on a graph. The set of all graphs whose $d$-separation relations are exactly the independence relations is the set of inferred graphs, our possible probabilistic models of the data. In general, this set will have more than one element, and graph with the same implied independence relations (i.e. $\triangle(G_1) = \triangle(G_2)$) are called a Markov equivalence class.

A crucial aspect is that we only used the independence relations between nodal variables for discovery. This work, [5], and [6], take first steps towards relaxing this assumption. More specifically, structural independence a 'structural' criterion similar to $d$-separation that characterizes independences in all distributions between *all* possible random variables, instead of just nodal random variables. In Section 9, we talk about the application of structural independence to structural causal models, a more specific model than a directed graph.

We are confident that this and further work will lead to better statistical methods for causal discovery and principled new probabilistic models. In Section 3, we talk more about the relationship to graphs.

Finally, we note that while [7] introduces deterministic nodes in a Bayesian network and a corresponding notion $D$-separation that captures structural independence for these models, $D$-separation leaves the functions that determine these nodes in general position. In the theory presented here, we fix a family of independent random elements and define structural independence for all random element that dependent only on this independent family. To further highlight this difference, our theory can be leveraged to define a $d$-separation criterion for arbitrary random variables defined on a graph, see [5], while $D$-separation does not do so.

# 3 Overview of the finite theory

We follow the exposition of [5] closely. Previous works also include [6] and [8]. The following definitions are valid throughout this section.

Let $I$ be a finite index set. Let $(\Omega, \mathcal{A})$ be a measurable space. Let $(U_i)_{i \in I}$ be a family of random elements on $\Omega$ with finite codomain. let

$$\triangle^{\times} := \Big\{ P : \mathcal{A} \to \mathbb{R} \mid P \text{ is a probability measure} \\ \text{and } (U_i)_{i \in I} \text{ is an independent family w.r.t. } P \Big\}$$

Note that here, in the finite case, it is not necessary to have a reference measure $\mathbb{P}$ s.t. $P$ is absolutely continuous w.r.t. $\mathbb{P}$. This is because we can choose a distribution $\mathbb{P}$ s.t. the pushforward $\mathbb{P}_U$ is the uniform distribution. Then all distributions on $\mathrm{Val}(U)$ are absolutely continuous w.r.t. $\mathbb{P}_U$.

We will now introduce the history and state the fundamental theorem without proofs. Proofs follow immediately from the general theory. Direct proofs can be read in [5], where everything is stated in the canonical space w.r.t. $U$ and without the use of measure theory.

First, we introduce generation, a sufficient condition for $X$ to be determined by $U_J$ given $C \in \sigma(U)$.

**Definition 3.1** (generation, history): Let $X$ be a random element and $C \in \sigma(U)$. We say that $J \subseteq I$ generates $X$ given $C$, if $\sigma(X|_C) \subseteq \sigma(U_J|_C)$ and $U(C) = U_J(C) \times U_{\overline{J}}(C)$.
The history of $X$ given $C$, is the subset-wise smallest $J \subseteq I$ that generates $X$ given $C$. The history of $X$ given $C$ is written $\mathcal{H}(X|C)$ and exists. For $A \in \mathcal{A}$, we also write $\mathcal{H}(A|C) \coloneqq \mathcal{H}(1_A|C)$.

The rectangle condition $U(C) = U_J(C) \times U_{\overline{J}}(C)$ is essential for the existence of the history and the correctness of the fundamental theorem. This condition reflects the requirement of independence of $U_J$ and $U_{\overline{J}}$ given $C$. Indeed, if this independence holds, the pushforward of $U$ under $P$ disintegrates into a product $P(U \in \cdot\,|C) = P(U_J \in \cdot \mid C) \times P(U_{\overline{J}} \in \cdot \mid C)$. If $P$ is chosen s.t. $P_U$ has no nontrivial nullsets, the support of $\mathbb{P}(U \in \cdot\,|C)$ is $U(C)$, while similarly, $\operatorname{supp} \mathbb{P}(U_J \in \cdot\,|C) = U_J(C)$ and $\operatorname{supp} \mathbb{P}(U_J \in \cdot\,|C) = U_{\overline{J}}(C)$. Now the mentioned product structure of $\mathbb{P}(U \in \cdot\,|C)$ implies that the support forms a Cartesian product. This exactly reflects $U(C) = U_J(C) \times U_{\overline{J}}(C)$. Conversely, it can be seen easily, that if $U(C) = U_J(C) \times U_{\overline{J}}(C)$, then $U_J$ and $U_{\overline{J}}$ are independent given $C$ in all distributions s.t. $U_J \perp\!\!\!\perp U_{\overline{J}}$.

**Theorem 3.2** (fundamental theorem): Let $X, Y$ and $Z$ be random elements with finite codomain. Then

$$\forall P \in \triangle^{\times} : X \perp\!\!\!\perp_P Y \mid Z \Leftrightarrow \forall z \in \operatorname{Val}(Z) : \mathcal{H}(X|\{Z = z\}) \cap \mathcal{H}(Y|\{Z = z\}) = \emptyset.$$

We can therefore define structural independence.

**Definition 3.3** (structural independence): Let $X, Y$ and $Z$ be random elements with finite codomain. Then $X$ and $Y$ are structurally independent given $Z$, if their conditional histories are disjoint. More precisely, we set

$$X \perp Y \mid Z :\Leftrightarrow \forall z \in \operatorname{Val}(Z) : \mathcal{H}(X|\{Z = z\}) \cap \mathcal{H}(X|\{Z = z\}) = \emptyset.$$

From this theory we can already express bayesnets. More precisely, in [5], from a directed acyclic graph $G$ we construct a family of background random variables $U$ and a set of random variables $(X_i)_{i \in I}$ that represent the nodes of $G$, s.t. the $d$-separation criterion on the graph is equivalent to structural independence of these random variables. Furthermore, $\mathcal{H}(X) \subseteq \mathcal{H}(Y)$ corresponds to $Y$ being an ancestor of $X$. In [5], this is called structural time. It is known that $d$-separation specifies a graph up to its skeleton (the undirected version of the graph), and certain arrows [2] Theorem 1.2.8. Therefore, $d$-separation and the ancestor relationship fully determines the graph. Here, we now have structural time and structural independence that generalize the ancestor relationship and $d$-separation respectively. In this sense, this theory also generalizes Pearl's theory.

# 4 Random index sets and random families

In this section let $I$ be an index set. Let $\mathfrak{P}(I)$ denote the powerset of $I$.

**Definition 4.1** (random index set): A random index set in $I$ on a measurable space $(\Omega, \mathcal{A})$ is a measurable mapping $\Omega \to \mathfrak{P}(I)$. Here, the powerset of $I$ is a measurable space endowed with the smallest $\sigma$-algebra that contains $\{K \subseteq I : i \in K\}$ for all $i \in I$.

**Lemma 4.2**: Let $J : \Omega \to \mathfrak{P}(I)$. Then $\sigma(J) = \sigma(\{i \in J\} : i \in I)$ and $J$ is measurable, i.e. a random index set, if and only if $\{i \in J\} \coloneqq \{\omega \in \Omega : i \in J(\omega)\}$ is measurable for all $i \in I$.

*Proof*: Since the $\sigma$-algebra on $\mathfrak{P}(I)$ (Definition 4.1) is generated by $\{K \subseteq I : i \in K\}$, we have

$$\begin{aligned}\sigma(J) &= \sigma(J^{-1}\{K \subseteq I, i \in K\} : i \in I) \\ &= \sigma(\{i \in J\} : i \in I).\end{aligned}$$

The second claim follows. □

**Definition 4.3**: Let $X = (X_i)_{i \in I}$ be a family of random elements. Let $J \subseteq I$. Then we denote by $X_J \coloneqq (X_i)_{i \in J}$, the restriction of this family to $J$. Canonically, $X_J$ is identified with the random element $(X_i)_{i \in J}(\omega) = (X_i(\omega))_{i \in J}$.

We now extend this standard notion to allow for random index sets. This means that the index set to which we restrict the family of random elements can vary with different outcomes $\omega \in \Omega$.

**Definition 4.4** (random family): Let $(X_i)_{i \in I}$ be a family of random elements on $\Omega$. Let $J : \Omega \to \mathfrak{P}(I)$ be a random index set. Then we define the random family $X_J$ by

$$\begin{aligned}X_J : \Omega &\to \bigcup_{K \subseteq I} \mathrm{Val}(X_K) \\ \omega &\mapsto X_{J(\omega)}(\omega) = (X_i(\omega))_{i \in J(\omega)}.\end{aligned}$$

We endow $\bigcup_{K \subseteq I} \mathrm{Val}(X_K)$ with the smallest $\sigma$-algebra that contains $\{x \in \mathrm{Val}(X_K) : K \subseteq I, i \in K \text{ and } x_i \in B\}$ for all $i \in I$ and measurable $B \subseteq \mathrm{Val}(X_i)$.

**Lemma 4.5**: The random family $X_J$ defined in Definition 4.4 is measurable and $\sigma(X_J) = \sigma(\{i \in J, X_i \in B\} : i \in I, B \subseteq \mathrm{Val}(X_i) \text{ measurable})$.

*Proof*: For $i \in I$ and $B \subseteq \mathrm{Val}(X_i)$ measurable, let

$$A_{i,B} = \{x \in \mathrm{Val}(X_K) : K \subseteq I, i \in K \text{ and } x_i \in B\} \subseteq \bigcup_{K \subseteq I} \mathrm{Val}(X_K)$$

By Definition 4.4, the sets $A_{i,B}$ generate the $\sigma$-algebra on $\mathrm{Val}(X_J)$. Therefore,

$$\begin{aligned}\sigma(X_J) &= \sigma\big(\quad X_J^{-1}(A_{i,B}) \quad \mid i \in I, B \subseteq \mathrm{Val}(X_i) \text{ measurable}\big) \\ &= \sigma(\{i \in J, X_i \in B\} \mid i \in I, B \subseteq \mathrm{Val}(X_i) \text{ measurable})\end{aligned}$$

In this expression, $J$ is a random index set and therefore $\{i \in J\}$ is measurable, and because $X_i$ is a random element and $B$ is measurable, $\{X_i \in B\}$ is measurable. Therefore, $\{i \in J, X_i \in B\}$ is measurable for any choice of $i$ and $B$. Therefore $X_J$ is measurable. □

**Corollary 4.6**: Let $X_J$ be the random family defined in Definition 4.4. Then $\sigma(X_J)$ is generated by the $\cap$-stable system

$$\left\{\{K \subseteq J, X_K \in B\} : K \subseteq I \text{ finite and } B = \bigtimes_{k \in K} B_k, \text{where } B_k \subseteq \mathrm{Val}(X_k) \text{ measurable}\right\}.$$

*Proof*: Clearly, this is system is $\cap$-stable, since

$$\begin{aligned}&\left\{K_1 \subseteq J, X_{K_1} \in B_1\right\} \cap \left\{K_2 \subseteq J, X_{K_2} \in B_2\right\} \\ &= \left\{K_1 \cup K_2 \subseteq J, X_{K_1 \cup K_2} \in \left(B_1 \times \mathrm{Val}\left(X_{K_2 \setminus K_1}\right)\right) \cap \left(B_2 \times \mathrm{Val}\left(X_{K_1 \setminus K_2}\right)\right)\right\}\end{aligned}$$

which is of the required form. It generates $\sigma(X_J)$ since it contains the generator given in Lemma 4.5. □

We now want to see how random families behave w.r.t. set operations on their random index sets.

**Notation 4.7**: For $J \subseteq I$ we denote its complement by $\overline{J} := I \setminus J$. When working with random index sets, set operations and relations are understood pointwise. Specifically, for random index sets $J$ and $K$, we define $\overline{J}(\omega) := \overline{J(\omega)}$, $(J \cap K)(\omega) := J(\omega) \cap K(\omega)$ and $(J \cup K)(\omega) := J(\omega) \cup K(\omega)$, as well as $J \subseteq K :\Leftrightarrow \forall \omega : J(\omega) \subseteq K(\omega)$.

**In the following**, let $(X_i)_{i \in I}$ be a family of random elements.

**Lemma 4.8**: Let $J$ be a random index set. Then $\sigma(J) \subseteq \sigma(X_J)$.

*Proof*: For $B = \mathrm{Val}(X_i)$, the set $\{\omega \in \Omega : i \in J(\omega) \text{ and } X_i(\omega) \in B\}$ is measurable and equal to $\{i \in J\}$. Sets of this form generate $\sigma(J)$ by Lemma 4.2 □

**Lemma 4.9**: Let $J$ and $K$ be random index sets, s.t. $J \subseteq K$. Then $\sigma(X_J) \subseteq \sigma(X_K, J)$.

*Proof*: Let $i \in I$ and $B \subseteq \mathrm{Val}(X_i)$ measurable. By Lemma 4.5, it suffices to show that $\{i \in J, X_i \in B\}$ is contained in $\sigma(X_K, J)$, since sets of this form generate $\sigma(X_J)$. Now, since $J \subseteq K$, we have

$$\{i \in J, X_i \in B\} = \underbrace{\{i \in J\}}_{\in \sigma(J)} \cap \underbrace{\{i \in K, X_i \in B\}}_{\in \sigma(X_K)} \in \sigma(X_K, J).$$

□

**Example 4.10**: In Lemma 4.9, it is essential to include the random index set $J$ in $\sigma(X_K, J)$. Let $I = \{1,2\}$, $\Omega = \{1,2,3\}^2$, and for $i \in I$ let $X_i = \pi_i : \Omega \to \{1,2,3\}$ be the projection on the $i$'th coordinate. Define random index sets $J$ and $K$ by

$$J(\omega) := \begin{cases} \{1\} & \text{if } \omega_2 = 1 \\ \emptyset & \text{else} \end{cases} \qquad K(\omega) := \begin{cases} \{1\} & \text{if } \omega_2 \in \{1,2\} \\ \emptyset & \text{else} \end{cases}$$

Clearly, $J \subseteq K$. Then by Lemma 4.8, $\{\omega \in \Omega : \omega_2 = 1\} = \{1 \in J\}$ is $\sigma(X_J)$-measurable. But it is not contained in $\sigma(X_K) = \sigma(\{\omega_1\} \times \{1,2\} : \omega_1 \in \{1,2,3\})$.

**Corollary 4.11**: Let $J$ and $K$ be random index sets. Then

$$\sigma(X_{J \cap K}) \subseteq \sigma(X_J, J \cap K) \cap \sigma(X_K, J \cap K)$$

*Proof*: The claim follows immediately from Lemma 4.9. □

**Example 4.12**: $\sigma(X_J) \cap \sigma(X_K) \subseteq \sigma(X_{J \cap K}, J, K)$ does not hold in general. Let $I = \{1,2\}, \Omega = \{1,2\}^2$ and for $i \in I$ let $X_i = \pi_i : \Omega \to \{1,2\}$ be the projection on the $i$'th coordinate. Define random index sets $J$ and $K$ by

$$J(\omega) := \begin{cases} \{1\} & \text{if } \omega_1 = \omega_2 \\ \{2\} & \text{else} \end{cases} \qquad K(\omega) := \begin{cases} \{2\} & \text{if } \omega_1 = \omega_2 \\ \{1\} & \text{else} \end{cases}$$

Then clearly, $\sigma(X_J) = \sigma(X_K) = \mathfrak{P}(\Omega)$, but since $J \cap K = \emptyset$, we have

$$\sigma(X_{J \cap K}, J, K) = \sigma(J, K) = \sigma(\{(1,1),(2,2)\}).$$

Therefore, $\sigma(X_J) \cap \sigma(X_K) \nsubseteq \sigma(X_{J \cap K}, J, K)$.

**Remark 4.13**: In Example 4.12, $\overline{J} = K$. Furthermore $\sigma(X_J) \cap \sigma(X_{\overline{J}}) = \mathfrak{P}(\Omega)$. Therefore we do not have $\sigma(X_J) \cap \sigma(X_{\overline{J}}) = \emptyset$ in general.

**Lemma 4.14**: Let $J$ and $K$ be random index sets. Then $\sigma(X_J, X_K) = \sigma(X_{J \cup K}, J, K)$.

*Proof*: '$\subseteq$' By symmetry it suffices to show that $\sigma(X_J) \subseteq \sigma(X_{J \cup K}, J, K)$. This follows directly from Lemma 4.9.

'⊇': Let $i \in I$ and $B \subseteq \mathrm{Val}(X_i)$ be measurable. By Lemma 4.8, we have $\sigma(J, K) \subseteq \sigma(X_J, X_K)$. Therefore, by Lemma 4.5, it suffices to show that $\{i \in J \cup K, X_i \in B\}$ is contained in $\sigma(X_J, X_K)$. Again, by Lemma 4.5,

$$\{i \in J \cup K, X_i \in B\} = \{i \in J, X_i \in B\} \cup \{i \in K, X_i \in B\} \in \sigma(X_J, X_K). \qquad \square$$

## 4.1 Almost sure relations

**In the following subsections**, let $\mathbb{P}$ be a probability measure on $(\Omega, \mathcal{A})$, and let $\mathcal{N} = \{A \in \mathcal{A} : \mathbb{P}(A) = 0\}$ denote the nullsets w.r.t. $\mathbb{P}$.

Recall that we say that an event $A \in \mathcal{A}$ holds almost surely, if $A^c \in \mathcal{N}$. For random variables $X, Y : \Omega \to \mathbb{R}$, we say that $X$ and $Y$ are almost surely equal and write $X \overset{\text{a.s.}}{=} Y$ if $\{X \neq Y\} \in \mathcal{N}$. Clearly, $\{X \neq Y\} = \{X - Y = 0\}$ is a measurable set, so this notion is well defined.

Similarly, we can define almost sure relations for events.

**Definition 4.1.1**: Let $A, B \in \mathcal{A}$ be events. We define the almost sure subset relation by $A \overset{\text{a.s.}}{\subseteq} B :\Leftrightarrow B \setminus A \in \mathcal{N}$ and almost sure equality by $A \overset{\text{a.s.}}{=} B :\Leftrightarrow A \overset{\text{a.s.}}{\subseteq} B \wedge B \overset{\text{a.s.}}{\subseteq} A$.

We can now use this notion to define almost sure relations for random index sets.

**Definition 4.1.2**: Let $J, K$ be random index sets. We define the almost sure subset relation by $J \overset{\text{a.s.}}{\subseteq} K :\Leftrightarrow \forall i \in I : \{i \in J\} \overset{\text{a.s.}}{\subseteq} \{i \in K\}$. and almost sure equality by $J \overset{\text{a.s.}}{=} K :\Leftrightarrow J \overset{\text{a.s.}}{\subseteq} K \wedge K \overset{\text{a.s.}}{\subseteq} J$.

Naively, we might want to define the almost sure relations for random index sets by $J \overset{\text{a.s.}}{\subseteq} K :\Leftrightarrow \{J \subseteq K\} \in \mathcal{N}$. But we run into two problems. First, it is not clear that $\{J \subseteq K\}$ is measurable, since we can now only write it as an uncountable intersection of measurable sets. $\{J \subseteq K\} = \bigcap_{i \in I} (\{i \notin J\} \cup \{i \in K\})$. Secondly, even if $\sigma(J) \subseteq \sigma(\mathcal{N})$, i.e. $J$ contains no nontrivial information, $\{J \neq \emptyset\}$ might be $\Omega$.

**Example 4.1.3**: Let $I = [0, 1]$ and $\Omega = [0, 1]^I$ endowed with the product $\sigma$-algebra. Let $\mathbb{P} = \bigtimes_{i \in I} \lambda|_{[0,1]}$ be the product of uniform distributions. Let $\mathcal{A}$ be the completed product $\sigma$-algebra on $\Omega$. Let $(U_i)_{i \in I} = (\pi_i)_{i \in I}$ be the coordinate projections. Define $J(\omega) = \{\omega_0\}$. Clearly, $J$ is measurable by Lemma 4.2, since for $i \in I$, $\{i \in J\} = \{\pi_0 = i\}$. Note that for any $i \in I$, $\{i \in J\}$ is a nullset, therefore $J \overset{\text{a.s.}}{=} \emptyset$ in the sense of Definition 4.1.2. But the set $\{J \neq \emptyset\} = \Omega$.

Example 4.1.3 shows that the naive definition of almost sure relations for random index sets is not useful for using random index sets to measure information content. We will later see that Definition 4.1.2 is the correct definition, because it is used successfully in the fundamental theorem (Theorem 6.4.1).

**Remark 4.1.4**: If $I$ is countable, the naive definition and Definition 4.1.2 coincide.

## 4.2 Minimality

We follow the setting of Section 4.1.

To show that the history in Section 5 is well defined, we need to show existence of certain (almost surely) minimal random index sets . In this section, investigate under which conditions a class of random index sets has an (almost surely) unique and minimal element.

**Terminology** Let $S$ be a set with a preorder $\leq$. Let $T \subseteq S$. We say $t \in T$ is minimal (w.r.t. $\leq$) in $T$, if $\forall s : s \leq t \Rightarrow t \leq s$. We say that there exists a (up to $\leq$) unique $t \in T$, if $\forall s \in T : s \leq t \wedge t \leq s$. Clearly, by taking a quotient w.r.t. the equivalence relation $s \sim t :\Leftrightarrow s \leq t \land t \leq s$, we recover a partial order on $S_{/\sim}$ where minimality and uniqueness correspond to the commonly known notions.

Therefore we can talk about the almost surely unique and minimal element in a class of events or random index sets. We say that a random index set [event] is almost surely minimal (unique) in a class of random index sets [events], if it is minimal (unique) w.r.t. the almost sure subset relation.

**Lemma 4.2.1**: Let $S$ be an index set. For $s \in S$, let $A_s \in \mathcal{A}$. Then there exists an almost surely unique and minimal set $B$, s.t. $\forall s \in S : A_s \overset{\text{a.s.}}{\subseteq} B$. Furthermore, there exists a countable set $S_0 \subseteq S$, s.t. $\bigcup_{s \in S_0} A_s \overset{\text{a.s.}}{=} B$.

*Proof*: The proof is trivial if $S$ is countable. We will directly construct $S_0 = \{s_n : n \in \mathbb{N}\}$ as a sequence and show that $\bigcup_{n \in \mathbb{N}} A_{s_n}$ can be chosen to be $B$. We inductively define $s_n$ by a sequence. For $n = 1$, let $s_n \in S$ be arbitrary. Let $n \in \mathbb{N}$ and $s_n$ be defined. Set $B_n := \bigcup_{m=1}^{n} A_{s_m}$ and $p_n := \sup_{s \in S} \mathbb{P}(A_s \setminus B_n)$. Choose $s_{n+1} \in S$, s.t. $\mathbb{P}\left(A_{s_{n+1}} \setminus B_n\right) \leq p_n - \frac{1}{n}$, We now prove the required properties.

1. $\forall s \in S : A_s \overset{\text{a.s.}}{\subseteq} B$: Let $s \in S$. We need to show $\mathbb{P}(A_s \setminus B) = 0$. Assume that $\mathbb{P}(A_s \setminus B) = p > 0$. Then $\mathbb{P}(A_s \setminus B_n) \geq \mathbb{P}(A_s \setminus B) = p$ because $B_n \subseteq B$. By the choice of $s_n$, $\mathbb{P}\left(A_{s_n} \setminus B_n\right) \geq \mathbb{P}(A_s \setminus B_n) - \frac{1}{n} \geq p - \frac{1}{n}$. But then with $B_0 = \emptyset$, $1 \geq \mathbb{P}(B) = \mathbb{P}\left(\bigcup_{n \in \mathbb{N}} A_{s_n} \setminus B_{n-1}\right) = \sum_{n \in \mathbb{N}} \mathbb{P}\left(A_{s_n} \setminus B_{n-1}\right) = \infty$, a contradiction.
2. Minimality and uniqueness: Let $C \in \mathcal{A}$, s.t. $\forall s \in S : A_s \overset{\text{a.s.}}{\subseteq} C$. Then clearly, $B = \bigcup_{n \in \mathbb{N}} A_{s_n} \overset{\text{a.s.}}{\subseteq} C$. Therefore, if $C$ is minimal too, i.e. $C \overset{\text{a.s.}}{\subseteq} B$, we have $B \overset{\text{a.s.}}{=} C$. $\square$

**Definition 4.2.2** (almost sure union and intersection): We call the set $B$ from Lemma 4.2.1 the almost sure union of the family $(A_s)_{s \in S}$ and write $\bigcup^{\text{a.s.}}_{s \in S} A_s := A$. We define the almost sure intersection of $(A_s)_{s \in S}$ by $\bigcap^{\text{a.s.}}_{s \in S} A_s := \left(\bigcup^{\text{a.s.}}_{s \in S} A_s^c\right)^c$. It then follows immediately that the almost sure intersection is the almost surely unique, maximal set that fulfills $\bigcap^{\text{a.s.}}_{s \in S} A_s \overset{\text{a.s.}}{\subseteq} A_s$ for all $s \in S$.

We now extend these notions to random index sets.

**Definition 4.2.3**: Let $S$ be an index set and $J_s$ be a random index set for all $s \in S$. We define the almost sure union of the family $(J_s)_{s \in S}$ to be the random index set $K$, s.t. $\{i \in K\} \overset{\text{a.s.}}{=} \bigcup^{\text{a.s.}}_{s \in S} \{i \in J_s\}$.

We denote $K$ by $\bigcup^{\text{a.s.}}_{s \in S} J_s$. Note that $K$ is only defined up to almost sure equality.

Similarly, we define the almost sure intersection of the family $(J_s)_{s \in S}$ and denote it by $\bigcap^{\text{a.s.}}_{s \in S} J_s$.

It is easy to see that the almost sure intersection is the almost surely unique, maximal random index set that fulfills $\bigcap^{\text{a.s.}}_{s \in S} J_s \overset{\text{a.s.}}{\subseteq} J_s$ for all $s \in S$. Similarly, the almost sure union is the almost surely unique, minimal random index set that fulfills $J_s \overset{\text{a.s.}}{\subseteq} \bigcup^{\text{a.s.}}_{s \in S} J_s$ for all $s \in S$.

With the almost sure intersection, we can now define $\{J \overset{\text{a.s.}}{=} K\}$, our substitute for the potentially non-measurable set $\{J = K\}$. Our definition will capture the notion that $J$ and $K$ contain the same information.

**Definition 4.2.4**: Let $J, K$ be random index sets. We define

$$\{J \overset{\text{a.s.}}{=} K\} := \bigcap_{i \in I}^{\text{a.s.}} (\{i \in J \cap K\} \cup \{i \in J^c \cap K^c\})$$

$$\{J \overset{\text{a.s.}}{\subseteq} K\} := \bigcap_{i \in I}^{\text{a.s.}} (\{i \notin J\} \cup \{i \in K\}).$$

Note that $\{J \overset{\text{a.s.}}{=} K\}$ and $\{J \overset{\text{a.s.}}{\subseteq} K\}$ are defined up to almost sure equality.

We can now apply Zorn's lemma to obtain minimal elements in certain classes of random index sets

**Lemma 4.2.5**: Let $\mathfrak{I}$ be a set of index functions s.t. for any totally ordered set $(S, \leq)$ and any random index set family $(J_s)_{s \in S}$ in $\mathfrak{I}$ with $\forall s, t \in S : s \leq t \Rightarrow J_s \overset{\text{a.s.}}{\subseteq} J_t$, we have $\bigcap_{s \in S}^{\text{a.s.}} J_s \in \mathfrak{I}$. Then there exists a almost surely subset wise minimal element in $\mathfrak{I}$. Furthermore, if $J, K \in \mathfrak{I}$ implies $J \cap K \in \mathfrak{I}$, this element is almost surely unique and is given by $\bigcap^{\text{a.s.}} \mathfrak{I}$.

*Proof*:

1. Existence: Apply Zorn's lemma to the random index sets modulo almost sure equality where the ordering is given by $J \leq K :\Leftrightarrow J \overset{\text{a.s.}}{\subseteq} K$.
2. Uniqueness: Assume that $J$ and $K$ are minimizers. If $J \overset{\text{a.s.}}{=} K$ does not hold, $J \cap K \overset{\text{a.s.}}{=} J$ does not hold and therefore $J$ is not a minimizer.
3. Let $J$ be the unique minimizer. Let $i \in I$. There exists a countable set $\mathfrak{I}_0 \subseteq \mathfrak{I}$, s.t. $\{i \in \bigcap^{\text{a.s.}} \mathfrak{I}_0\} \overset{\text{a.s.}}{=} \{i \in \bigcap^{\text{a.s.}} \mathfrak{I}\}$. Enumerate $\mathfrak{I}_0$ by $J_n$ and set $K_n := J_1 \cap \ldots \cap J_n$. Note that $K_n \in \mathfrak{I}$ by assumption and $n \geq m \Rightarrow K_m \subseteq K_n$. Then $K := \bigcap^{\text{a.s.}} \mathfrak{I}_0 = \bigcap_{n \in \mathbb{N}} K_n \in \mathfrak{I}$ by assumption (where we use the dual order on $\mathbb{N}$). Now $\{i \in J\} \overset{\text{a.s.}}{\subseteq} \{i \in K\} = \{i \in \bigcap^{\text{a.s.}} \mathfrak{I}\}$ holds by construction. Furthermore, $\{i \in J\} \overset{\text{a.s.}}{\supseteq} \{i \in \bigcap^{\text{a.s.}} \mathfrak{I}\}$ holds by definition of $\bigcap^{\text{a.s.}} \mathfrak{I}$ and the fact that $J \in \mathfrak{I}$. Since $i \in I$ was arbitrary, we have $J \overset{\text{a.s.}}{=} \bigcap^{\text{a.s.}} \mathfrak{I}$. □

## 4.3 Generated $\sigma$-algebras

**In the following**, let $\mathbb{P}$ be a probability measure on $(\Omega, \mathcal{A})$ and $\mathcal{N} = \{A \in \mathcal{A} : \mathbb{P}(A) = 0\}$ the nullsets w.r.t. $\mathbb{P}$. Let $\mathcal{A}$ be complete w.r.t. $\mathbb{P}$, i.e. $\mathcal{N} \subseteq \mathcal{A}$.

We now want to see how the almost sure relations for random index sets behave w.r.t. the generated $\sigma$-algebras of both themselves and their corresponding random families.

**Notation 4.3.1** (completed $\sigma$-algebras): For ease of notation, we will assume all $\sigma$-algebras to be complete w.r.t. the reference measure $\mathbb{P}$. We already assume that $\mathcal{A}$ is complete and whenever we write $\sigma_{\mathcal{N}}(\cdot)$, we mean the generated $\sigma$-algebra, completed w.r.t. $\mathbb{P}$, i.e. $\sigma_{\mathcal{N}}(\cdot) = \sigma(\cdot, \mathcal{N})$. This is a very well behaved operation, since for any set systems $\mathcal{B}_1, \mathcal{B}_2$, we have $\sigma(\mathcal{B}_1) \subseteq \sigma(\mathcal{B}_2) \Rightarrow \sigma_{\mathcal{N}}(\mathcal{B}_1) \subseteq \sigma_{\mathcal{N}}(\mathcal{B}_2)$. Therefore all results proven in Section 4 about generated $\sigma$-algebras carry over to $\sigma_{\mathcal{N}}(\cdot)$.

**In the remainder of this work**, we write $\sigma(\cdot)$ for $\sigma_{\mathcal{N}}(\cdot)$ for ease of notation.

**Remark 4.3.2**: Similar to the space of $L^p$ functions, we can take the quotient w.r.t. the $\sigma$-ideal $\mathcal{N}$ in the space of sub $\sigma$-algebras of $\mathcal{A}$, events in $\mathcal{A}$, but this is not strictly necessary, and we prefer to work with the completed $\sigma$-algebras. We will instead keep making use of almost sure relations by writing a.s. on top of the relations. We can also take the quotient in the space of random index sets w.r.t. almost sure equality, which allows us to talk about a minimal and unique random index set (equivalence class) in Lemma 4.2.5. For clarity we will still denote this by almost sure minimality and uniqueness.

The following lemmas prove that random index sets and random elements behave as expected w.r.t. almost sure relations.

This first lemma shows that almost surely equal random index sets contain the same information.

**Lemma 4.3.3**: Let $J$ and $K$ be random index sets, s.t. $J \overset{\text{a.s.}}{=} K$, then $\sigma(J) \overset{\text{a.s.}}{=} \sigma(K)$. Moreover, if $J \overset{\text{a.s.}}{=} \emptyset$, then $\sigma(J) = \sigma(\emptyset)$.

*Proof*: By definition, $\{i \in J\} \overset{\text{a.s.}}{=} \{i \in K\}$ for all $i \in I$. These sets and $\mathcal{N}$ generate $\sigma(J)$ and $\sigma(K)$ respectively. If $J \overset{\text{a.s.}}{=} \emptyset$, then $\{i \in J\} \in \mathcal{N}$. $\square$

This lemma shows that a random index set that contains no information produces a random family that contains no information.

**Lemma 4.3.4**: Let $J$ be a random index set s.t. $J \overset{\text{a.s.}}{=} \emptyset$. $\sigma(U_J) = \sigma(\emptyset)$.

*Proof*: Let $i \in I$ and $B \subseteq \mathrm{Val}(X_i)$ measurable. It suffices to show that $\mathbb{P}\{i \in J, X_i \in B\} = 0$, since sets of this form and the nullsets $\mathcal{N}$ generate $\sigma(U_J)$ and therefore $\sigma(U_J) \subseteq \sigma(\mathcal{N}) = \sigma(\emptyset)$. This follows from $J \overset{\text{a.s.}}{=} \emptyset \Rightarrow \forall i \in I : \{i \in J\} \in \mathcal{N}$. $\square$

**Lemma 4.3.5**: Let $J$ and $K$ be random index sets s.t. $J \overset{\text{a.s.}}{=} K$. Then $\sigma(U_J) = \sigma(U_K)$.

*Proof*: Clearly, $J \cap K \overset{\text{a.s.}}{=} J$. We can therefore w.l.o.g. assume $J \subseteq K$. By Lemma 4.14, $\sigma(U_J, U_{K \setminus J}) = \sigma(U_K, J, K \setminus J)$. By Lemma 4.3.4, $\sigma(U_{K \setminus J}) = \sigma(\emptyset)$. By Lemma 4.3.3, $\sigma(J) = \sigma(K)$ and $\sigma(K \setminus J) = \sigma(\emptyset)$. Therefore, $\sigma(U_J) = \sigma(U_K, K) = \sigma(U_K)$ Lemma 4.8. $\square$

In most of the lemmas about relationships between $\sigma$-algebras of random families, we had to include the random index set in the $\sigma$-algebra (c.f. Lemma 4.9). Once we introduce a random element w.r.t. which the involved random index sets are measurable, this need disappears.

**Lemma 4.3.6**: Let $S$ be an index set. Let $Z$ be a random element. Let $(J_s)_{s \in S}$ be a family of $\sigma(Z)$-measurable random index sets. Let $J = \bigcup^{\text{a.s.}}_{s \in S} J_s$. Then $\sigma(U_J, Z) = \sigma(U_{J_s}, Z : s \in S)$.

*Proof*: '$\supseteq$': Follows immediately by Lemma 4.9 and the $\sigma(Z)$-measurability of all involved random index sets.
'$\subseteq$': Let $i \in I$ and $B \subseteq \mathrm{Val}(X_i)$ measurable. It suffices to show that $\{i \in J, X_i \in B\} \in \sigma\left(U_{J_s}, Z : s \in S\right)$, since these sets generate $\sigma(U_J)$. By Definition 4.2.4 and Lemma 4.2.1, there is a countable subset $S_0 \subseteq S$, s.t. $\{i \in J\} \overset{\text{a.s.}}{=} \bigcup_{s \in S_0} \{i \in J_s\}$. Now clearly, $\{i \in J, X_i \in B\} \overset{\text{a.s.}}{=} \bigcup_{s \in S_0} \{i \in J_s, X_i \in B\} \in \sigma\left(U_{J_s}, Z : s \in S\right)$. $\square$

In Lemma 4.3.6, $Z$ will naturally appear as the random element on which we condition.

# 5 Construction of the history

Let $(\Omega, \mathcal{A}, \mathbb{P})$ be a complete probability space, and $\mathcal{N}$ denote the nullsets w.r.t. $\mathbb{P}$. Recall that $\sigma(X) = \sigma_{\mathcal{N}}(X)$ denotes the completion w.r.t. $\mathbb{P}$ of the generated $\sigma$-algebra of $X$ (Notation 4.3.1).

Let $I$ be an arbitrary index set. Let $(U_i)_{i \in I}$ be a family of random elements. Define

$$\begin{aligned} \triangle^{\times} := \Big\{ P : \mathcal{A} \to \mathbb{R} \mid\ & P \text{ is a probability measure,} \\ & (U_i)_{i \in I} \text{ is an independent family w.r.t. } P \\ & \text{and } P \sim \mathbb{P} \Big\} \end{aligned}$$

**In the following** $X, Y$ and $Z$ will be arbitrary random elements.

We can now use random index sets to describe the dependence on $U$ as $Z$ varies. Let $\Sigma$ denote the set of all (complete) sub $\sigma$-algebras of $\mathcal{A}$. Our goal is to construct a map we will call history, $\mathcal{H}(\cdot \mid \cdot) : \Sigma \times \Sigma \to \mathfrak{P}(I)^\Omega$ that measures the dependence on $U$ of a $\sigma$-algebra given another $\sigma$-algebra with a random index set. We write $\mathcal{H}(X|Z)$ for $\mathcal{H}(\sigma(X)|\sigma(Z))$.

The name history is motivated by the fact that to determine the outcome of $X$, we need to know the outcomes of $U_i$ for all $i \in \mathcal{H}(X|\emptyset)$. Morally, $U_i$ need to 'have happened' for $X$ to be determined.

The conditional history $\mathcal{H}(X|Z)$ is harder to interpret. It characterizes the dependence of $X$ on $U$ given $Z$. Morally, given that $Z$ has happened, $X$ depends only on $U_i$ for $i \in \mathcal{H}(X|Z)$, and is independent of $U_i$ for $i \notin \mathcal{H}(X|Z)$.

If $\mathcal{H}(X|Z) = \emptyset$, then $X$ does not depend on anything given $Z$, i.e. $\sigma(X) \subseteq \sigma(Z)$. All the information that $X$ provides is already contained in $Z$. For polish spaces, this is equivalent to $X$ being a deterministic function of $Y$, i.e. there is a measurable function $f : \mathrm{Val}(Y) \to \mathrm{Val}(X)$, s.t. $X = f(Y)$.

**Desiderata 5.1** (Desiderata for the history):

1. $\mathcal{H}(X|Z)$ is a $\sigma(Z)$-measurable random index set.
   Morally, the set $\{i \in \mathcal{H}(X|Z)\}$ represent the event that $X$ depends (at least) on $U_i$ given $Z$. Since we are conditioning on $Z$, we require that $\{i \in \mathcal{H}(X|Z)\}$ is $\sigma(Z)$-measurable. This implies that $\mathcal{H}(X|Z)$ needs to be a $\sigma(Z)$ random index set.
2. Characterizes independence: $\forall P \in \triangle^\times : X \perp\!\!\!\perp_P Y \mid Z \Leftrightarrow \mathcal{H}(X|Z) \cap \mathcal{H}(Y|Z) \overset{\text{a.s.}}{=} \emptyset$.
   We argue for both directions separately. '$\Leftarrow$' means that whenever $X$ and $Y$ depend on disjoint sets of $U_i$ given $Z$, they are independent given $Z$, regardless of the specific distribution $P \in \triangle^\times$. This is because we used the assumption in the definition of $\triangle^\times$ that $U_i$ are independent w.r.t. $P$.
   '$\Rightarrow$' says that disjointness of the histories is the weakest condition that implies independence given $Z$ regardless of distribution. Morally, it means that the only independence (regardless of distribution) we can achieve is the one that is provided by the assumption that $U_i$ are independent and an application of a conditional version of the standard independence theorem that functions of disjoint subfamilies of an independent family are independent.
3. Monotonicity: If $\sigma(X) \subseteq \sigma(Y)$ then $\mathcal{H}(X|Z) \overset{\text{a.s.}}{\subseteq} \mathcal{H}(Y|Z)$.
   Morally, gaining information content can only increase the set of $U_i$ on which $X$ depends given $Z$.
4. $\mathcal{H}(U|Z) \cap J \overset{\text{a.s.}}{\subseteq} \mathcal{H}(U_J|Z)$ for all random index sets $J$.
   Morally, the part of $J$ that is relevant for $U$, is also relevant for $U_J$.

Note that condition 4. anchors $i$ to $U_i$. Otherwise, we could take any permutation $\tau$ of $I$, and define a new history by $\{i \in \mathcal{H}'(X|Z)\} := \{\tau(i) \in \mathcal{H}(X|Z)\}$ that satisfies desiderata 1-3.

In Theorem 7.16, we will see that our construction is the (pointwise) maximal map that fulfills Desiderata 5.1. Maximality is a desirable property, since the disjointness of $\mathcal{H}(X|Z)$ and $\mathcal{H}(Y|Z)$ in Desiderata 5.1.2 becomes harder to fulfill if $\mathcal{H}(X|Z)$ and $\mathcal{H}(Y|Z)$ are larger.

Our construction will mimic the finite case (Definition 3.1). Specifically, we first define what it means for a random index set $J$ to be sufficient to determine for $X$ given $Z$. Then we show that there exists an almost surely minimal and unique sufficient $J$.

**Definition 5.2** (disintegration, generation): Let $J$ be a $\sigma(Z)$-measurable random index set. We say that $J$ disintegrates $Z$, if $\forall P \in \triangle^{\times} : U_J \perp\!\!\!\perp_P U_{\overline{J}} \mid Z$.
We say that $J$ generates $X$ given $Z$, if $\sigma(X) \subseteq \sigma(U_J, Z)$ and $J$ disintegrates $Z$.

We can see that disintegration is actually a well-behaved notion and we don't need to quantify over all distributions. Rather, it suffices to check the condition for any distribution, in particular we can choose the reference measure $\mathbb{P}$.

**Lemma 5.3**: $J$ disintegrates $Z$ if and only if $U_J \perp\!\!\!\perp_{\mathbb{P}} U_{\overline{J}} \mid Z$.

*Proof*: '$\Rightarrow$' Follows from the definition of disintegration.
'$\Leftarrow$': Let $P \in \triangle^{\times}$. By Lemma A.8 there is a density $\varphi$, s.t. $P = \varphi \cdot \mathbb{P}$ and a family of densities $(\varphi_n)_{n \in \mathbb{N}_0}$ and a family of indices $(i_n)_{n \in \mathbb{N}}$ s.t. $\mathbb{E}(\varphi_0 | U) = 1$, and $\varphi_n$ is $U_{i_n}$-measurable. Furthermore, $\prod_{n \in \mathbb{N}_0} \varphi_n$ converges $\mathbb{P}$-a.s. to $\varphi$. Let $J' = J \cap \{i_n : n \in \mathbb{N}\}$ and $\overline{J}' = \overline{J} \cap \{i_n : n \in \mathbb{N}\}$ We define $\varphi_J = \prod_{i \in J'} \varphi_i$ and $\varphi_{\overline{J}} = \prod_{i \in \overline{J}'} \varphi_i$. Since $\varphi_J = \prod_{n \in \mathbb{N}} \left(1_{\{i_n \in J\}} \varphi_{i_n} + 1_{\{i_n \notin J\}}\right)$, $\varphi_J$ is $\sigma(U_J)$-measurable and $\varphi_{\overline{J}}$ is $\sigma(U_{\overline{J}})$-measurable. Let $\mathbb{E}$ be the expectation w.r.t. $\mathbb{P}$. Let $A \in \sigma(U_J)$ and $B \in \sigma(U_{\overline{J}})$. Then $P(A|Z)P(B|Z) = P(A, B|Z) \Leftrightarrow \mathbb{E}(\varphi 1_A | Z)\mathbb{E}(\varphi 1_B | Z) = \mathbb{E}(\varphi | Z)\mathbb{E}(\varphi 1_A 1_B | Z)$. The claim follows by $U_J \perp\!\!\!\perp_{\mathbb{P}} U_{\overline{J}} \mid Z$ through the equalities $\mathbb{E}(\varphi 1_A | Z) = \mathbb{E}(\varphi_J 1_A | Z)\mathbb{E}(\varphi_{\overline{J}})$, $\mathbb{E}(\varphi 1_B | Z) = \mathbb{E}(\varphi_J | Z)\mathbb{E}(\varphi_{\overline{J}} 1_B | Z)$ and $\mathbb{E}(\varphi 1_A 1_B | Z) = \mathbb{E}(\varphi_J 1_A | Z)\mathbb{E}(\varphi_{\overline{J}} 1_B | Z)$. □

The following corollary shows that generation behaves well with random index sets defined up to almost sure equality. Therefore, it makes sense to ask for the almost surely minimal generating random index set.

**Corollary 5.4**: Let $J$ and $K$ be random index sets s.t. $J \overset{\text{a.s.}}{=} K$. Then $J$ generates $X \mid Z$ if and only if $K$ generates $X \mid Z$.

*Proof*: Follows immediately from Lemma 4.3.5, because Definition 5.2 only uses $J$ through $\sigma(U_J)$, $\sigma(U_{\overline{J}})$. □

To prove that a minimal generating index set exists, we want to use Lemma 4.2.5, the application of Zorn's lemma to random index sets. For this we need to show that the class of random index sets that generate $X$ given $Z$ is closed under intersections and under taking intersections of chains.

First, we look at the intersection.

**Lemma 5.5**: Let $J$ and $K$ be random index sets that disintegrate $Z$. Then $J \cap K$ disintegrate $Z$.

*Proof*: We let $(J_1, J_2, J_3, J_4) \coloneqq (J \cap K, J \setminus K, K \setminus J, K)$ and $J_{ij} \coloneqq J_i \cup J_j$. By assumption, $J_{12} = J$ and $J_{13} = K$ disintegrate $Z$. We need to show that $J_1 = J \cap K$ disintegrates $Z$.

Let $P \in \triangle^{\times}$. It suffices to show that the family $(U_{J_i})_{i=1}^4$ is independent given $Z$. This is equivalent to the conditional independence of $(U_{J_i}, Z)_{i=1}^4$. By Lemma 4.14 and the $\sigma(Z)$-measurability of $J$ and $K$, we have $\sigma\left(U_{J_1}, U_{J_2}, Z\right) = \sigma\left(U_{J_{12}}, Z\right)$ and likewise for other choices of indices. By disintegration w.r.t $J_{12}$ and $J_{13}$, we have

$$\left(U_{J_{12}}, Z\right) \perp\!\!\!\perp_P \left(U_{J_{34}}, Z\right) \mid Z \tag{1}$$

$$\left(U_{J_{13}}, Z\right) \perp\!\!\!\perp_P \left(U_{J_{24}}, Z\right) \mid Z \tag{2}$$

For $i \in \{1, ..., 4\}$ let $A_i \in \sigma\left(U_{J_i}, Z\right)$. It suffices to show $A_1 \perp\!\!\!\perp_P A_2 \cap A_3 \cap A_4 \mid Z$, since sets of this form are a $\cap$-stable system that generates $\sigma(U_{234}, Z)$.

Now clearly, $P\left(\bigcap_{i=1}^{4} A_i | Z\right) \overset{(1)}{=} P(A_1 \cap A_2 \mid Z) P(A_3 \cap A_4 \mid Z) \overset{(2)}{=} \prod_{i=1}^{4} P(A_i | Z)$. □

**Corollary 5.6**: Let $J$ and $K$ be random index sets that disintegrate $Z$. Then $J \cup K$ disintegrate $Z$.

*Proof*: By symmetry of disintegration, we have $J^c$ and $K^c$ disintegrate $Z$. By Lemma 5.5, $J^c \cap K^c$ disintegrates $Z$. Again by the symmetry of disintegration, we have $J \cup K = (J^c \cap K^c)^c$ disintegrates $Z$. □

**Lemma 5.7**: Let $J$ and $K$ generate $X$ given $Z$. Then $J \cap K$ generates $X$ given $Z$.

*Proof*: By Lemma 5.5, $J \cap K$ disintegrates $X$ given $Z$. By the definition of generation, we have $\sigma(X) \subseteq \sigma(U_J, Z) \cap \sigma(U_K, Z)$. It suffices to show that the right hand side is a subset of $\sigma(U_{J \cap K}, Z)$. Let $A \in \sigma(U_J, Z) \cap \sigma(U_K, Z)$.

Since $U_{J \cap K} \perp\!\!\!\perp_{\mathbb{P}} U_{J \setminus K} \mid Z$, we have $\mathbb{P}(U_{J \setminus K} \in \cdot \mid Z) = \mathbb{P}(U_{J \setminus K} \in \cdot \mid U_{J \cap K}, Z)$. Since the same holds true for $K \setminus J$, we get $U_{J \setminus K} \perp\!\!\!\perp_{\mathbb{P}} U_{K \setminus J} \mid U_{J \cap K}, Z$ and therefore $(U_{J \setminus K}, U_{J \cap K}, Z) \perp\!\!\!\perp_{\mathbb{P}} (U_K, U_{J \cap K}, Z) \mid U_{J \cap K}, Z$. By Lemma 4.14, we have $\sigma(U_J, Z) \perp\!\!\!\perp_{\mathbb{P}} \sigma(U_K, Z) \mid U_{J \cap K}, Z$ Since $A \in \sigma(U_J, Z) \cap \sigma(U_K, Z)$, we have $A \perp\!\!\!\perp_{\mathbb{P}} A \mid U_{J \cap K}, Z$. This implies $A \in \sigma(U_{J \cap K}, Z)$. □

It remains to show that generation is closed under chains, i.e. that a descending chain of generating random index sets that generate $X$ given $Z$ can be intersected to obtain a generating random index set.

**Lemma 5.8**: Let $(S, \leq)$ be a totally ordered set. Let $(J_s)_{s \in S}$ be a family of random index sets that disintegrate $Z$. s.t. $\forall s, t \in S : s \leq t \Rightarrow J_s \overset{\text{a.s.}}{\subseteq} J_t$. Then $\bigcap^{\text{a.s.}}_{s \in S} J_s$ disintegrates $Z$.

*Proof*: Let $J := \bigcap^{\text{a.s.}}_{s \in S} J_s$. Let $P \in \triangle^{\times}$. We need to show that $U_J \perp\!\!\!\perp_P U_{\overline{J}} \mid Z$. For this let $A \in \sigma(U_J, Z)$ and $s \in S$ and $B \in \sigma\left(U_{\overline{J_s}}, Z\right)$. Since $\overline{J} = \bigcup^{\text{a.s.}}_{s \in S} \overline{J_s}$ and $\bigcup_{s \in S} \sigma\left(U_{\overline{J_s}}, Z\right)$ is a $\cap$-stable system that generates $\sigma(U_{\overline{J}}, Z)$ (Lemma 4.3.6), it suffices to show $A \perp\!\!\!\perp_P B \mid Z$ for all $B \in \sigma\left(U_{\overline{J_s}}, Z\right)$. Now clearly, $A \in \sigma\left(U_{J_s}, Z\right)$, so that $A \perp\!\!\!\perp_P B \mid Z$ because $J_s$ disintegrates $Z$. □

**Lemma 5.9**: Let $(S, \leq)$ be a totally ordered set. Let $(J_s)_{s \in S}$ be a family of random index sets that generate $X$ given $Z$. s.t. $\forall s, t \in S : s \leq t \Rightarrow J_s \overset{\text{a.s.}}{\subseteq} J_t$. Then $\bigcap^{\text{a.s.}}_{s \in S} J$ generates $Z$ given $X$

*Proof*: Let $J := \bigcap^{\text{a.s.}}_{s \in S} J_s$. By Lemma 5.8, $J$ disintegrates $X$ given $Z$, therefore $U_J \perp\!\!\!\perp_{\mathbb{P}} U_{\overline{J}} \mid Z$. By the definition of generation, we have $\sigma(X) \subseteq \bigcap_{s \in S} \sigma\left(U_{J_s}, Z\right)$. It suffices to show that $\bigcap_{s \in S} \sigma\left(U_{J_s}, Z\right) \overset{\text{a.s.}}{\subseteq} \sigma(U_J, Z)$.

Let $s \in S$, then by Lemma 4.3.6, $\sigma\left(U_{\overline{J_s}}, Z\right) \subseteq \sigma(U_{\overline{J}}, Z)$ and $\sigma(U_{\overline{J}}, Z) \subseteq \sigma\left(U_{J_s \setminus J}\right)$. Therefore $U_J \perp\!\!\!\perp_{\mathbb{P}} U_{\overline{J_s}} \mid Z$, and $U_J \perp\!\!\!\perp_{\mathbb{P}} U_{J_s \setminus J} \mid Z$.
This implies $\mathbb{P}\left(U_{\overline{J_s}} \in \cdot \mid Z\right) = \mathbb{P}\left(U_{\overline{J_s}} \in \cdot \mid Z, U_J\right)$ and $\mathbb{P}\left(U_{J_s \setminus J} \in \cdot \mid Z\right) = \mathbb{P}\left(U_{J_s \setminus J} \in \cdot \mid Z, U_J\right)$.
Since $J_s$ disintegrates $Z$ and by Lemma 4.3.6, we have $U_{J_s \setminus J} \perp\!\!\!\perp_{\mathbb{P}} U_{\overline{J_s}} \mid Z, U_J$.
Therefore $U_{J_s} \perp\!\!\!\perp_{\mathbb{P}} U_{\overline{J_s} \cup J} \mid Z, U_J$. Finally, let $A \in \bigcap_{s \in S} \sigma\left(U_{J_s}, Z\right)$. By the previous independence statement, $A \perp\!\!\!\perp_{\mathbb{P}} U_{\overline{J_s} \cup J} \mid Z, U_J$ for all $s \in S$. Since $\bigcup_{s \in S} \sigma\left(U_{\overline{J_s} \cup J}\right)$ is a $\cap$-stable generator of $\sigma(U_{\overline{J} \cup J}) = \sigma(U)$, we get $A \perp\!\!\!\perp_P A \mid Z, U_J$, and therefore $A \in \sigma(U_J, Z)$. □

**Theorem 5.10**: There exists an almost surely unique minimal generating random index set of $X$ given $Y$, that is given by $\bigcap^{\text{a.s.}} \{J : \Omega \to \mathfrak{P}(I) \mid J \text{ generates } X \text{ given Z}\}$.

*Proof*: Apply Zorn's lemma (Lemma 4.2.5) with the closure of generation under intersection (Lemma 5.7) and chains (Lemma 5.9) . □

**Definition 5.11** (history): The almost surely unique minimal generating random index set of $X$ given $Y$ is called the history of $X$ given $Y$. We write

$$\mathcal{H}(X|Z) := \overset{\text{a.s.}}{\bigcap} \{J : \Omega \to \mathfrak{P}(I) \mid J \text{ generates } X \text{ given Z}\}$$

Theorem 5.10 states that $\mathcal{H}(X|Z)$ generates $X$ given $Z$.

We are now ready to define structural independence in terms of histories.

**Definition 5.12**: $X$ and $Y$ are structurally independent given $Z$, if their histories are almost surely disjoint. More precisely,

$$X \perp Y \mid Z :\Leftrightarrow \mathcal{H}(X|Z) \cap \mathcal{H}(Y|Z) \overset{\text{a.s.}}{=} \emptyset.$$

# 6 The fundamental theorem of structural independence

The goal of this section is to prove that structural independence (Definition 5.12) characterizes independence in all product distributions $P \in \triangle^\times$ (Desiderata 5.1.2). The proof simplifies in the finite case and was presented in [5]. In the finite case, it is also possible to prove the statement using polynomials and factoring lemmas, see [6]. This is no longer possible in the infinite case, so we need to use a more direct approach.

Our goal is to prove

$$\begin{aligned} \forall P \in \triangle^\times : X \perp\!\!\!\perp_P Y \mid Z \Leftrightarrow X \perp Y \mid Z \\ :\Leftrightarrow \mathcal{H}(X|Z) \cap \mathcal{H}(Y|Z) \overset{\text{a.s.}}{=} \emptyset. \end{aligned}$$

This is Desiderata 5.1.2 for the history map. One direction is essentially implied directly by the definition of generation.

**Theorem 6.1** (Soundness of structural independence): Structural independence is sound: For any possible realization of a product probability distribution $P \in \triangle^\times$, structural independence implies independence. More precisely, $X \perp Y \mid Z \Rightarrow \forall P \in \triangle^\times : X \perp\!\!\!\perp_P Y \mid Z$.

*Proof*: Let $P \in \triangle^\times$. Since $J := \mathcal{H}(X|Z)$ generates $X$ given $Z$, we have $\sigma(X) \subseteq \sigma(U_J, Z)$ and $U_J \perp\!\!\!\perp_P U_{\overline{J}} \mid Z$. Using $\mathcal{H}(Y|Z) \overset{\text{a.s.}}{\subseteq} \overline{J}$ and applying Lemma 4.3.6, we get $\sigma(Y) \subseteq \sigma\left(U_{\mathcal{H}(Y|Z)}, Z\right) \subseteq \sigma(U_{\overline{J}}, Z)$. The claim now follows by the independence of $U_J$ and $U_{\overline{J}}$ given $Z$. □

To prove the other direction, completeness, we define a dual notion to the history and generation, irrelevance. While history and generation talk about depending, irrelevance talks about there being no information. To motivate this notion we look at a simple lemma.

**Lemma 6.2**: Let $i \in I$ and $P, Q \in \triangle^\times$, s.t. there is a $\sigma(U_i)$-measurable positive density $f : \Omega \to \mathbb{R}_{>0}$, s.t. $P = f \cdot Q$. Then for all $A \in \sigma(X)$ we have $P(A|Z)(\omega) = Q(A|Z)(\omega)$ for a.e. $\omega \in \left\{i \in \overline{\mathcal{H}(X|Z)}\right\}$.

*Proof*: Let $J := \mathcal{H}(X|Z)$ It suffices to show that $P(A|Z) \overset{\text{a.s.}}{=} Q(A|Z)$ for all $A \in \sigma(U_J, Z)|_{\{i \in \overline{J}\}}$. Let $A = B \cap C$, where $B \in \sigma(U_J)$ and $C \in \sigma(Z)|_{\{i \in \overline{J}\}}$. Since sets of this form are $\cap$-stable and generate $\sigma(U_J, Z)|_{\{i \in \overline{J}\}}$, and $\{A \in \mathcal{A} : P(A|Z) \overset{\text{a.s.}}{=} Q(A|Z)\}$ is a Dynkin system, it suffices to show the statement for such $A$. Note that $1_C f \perp\!\!\!\perp_P B \mid Z$, since $J$ disintegrates $Z$ and $\sigma(1_C f) \subseteq \sigma(U_{\overline{J}})$ and $B \in \sigma(U_J)$. Let $E$ denote the expectation w.r.t. $Q$. Now $P(A|Z) \overset{\text{a.s.}}{=} \frac{E(1_C f 1_B|Z)}{E(f|Z)} \overset{\text{a.s.}}{=} \frac{E(1_C f|Z) E(1_B|Z)}{E(f|Z)} \overset{\text{a.s.}}{=} 1_C E(1_B|Z) \overset{\text{a.s.}}{=} Q(A|Z)$. □

This lemma tells us that when we change the distribution of $U_i$, we cannot change $P(A|Z)$ in the region $\{i \notin \mathcal{H}(X|Z)\}$. This is a dual notion to the dependence of history. This motivates the following definitions.

## 6.1 The random index set of irrelevance

**Definition 6.1.1**: Let $K \subseteq I$. Let

$$\triangle_K^{\times 2} \coloneqq \{(P, Q) \in \triangle^\times \times \triangle^\times : \exists f : \Omega \to \mathbb{R}_{>0}, \sigma(U_K)\text{-measurable, s.t. } Q = f \cdot P\}$$

For $i \in I$, set $\triangle_i^{\times 2} = \triangle_{\{i\}}^{\times 2}$.

**Definition 6.1.2**: Let $J$ be a $\sigma(Z)$-measurable random index set. We say $J$ is irrelevant to $X$ given $Z$, if for all $i \in I$, $(P, Q) \in \triangle_i^{\times 2}$ and $A \in \sigma(X)$, we have $P(A|Z)(\omega) = Q(A|Z)(\omega)$ for a.e. $\omega \in \{i \in J\}$.

Of course, irrelevance extends to arbitrary subsets of $I$.

**Lemma 6.1.3**: Let $J$ be irrelevant to $X$ given $Z$. Let $K \subseteq I$. Then for $(P, Q) \in \triangle_K^{\times 2}$ and $A \in \sigma(X)$, we have $P(A|Z)(\omega) = Q(A|Z)(\omega)$ for a.e. $\omega \in \{K \overset{\text{a.s.}}{\subseteq} J\}$.

*Proof*: Let $\varphi$ be the probability density of $Q$ w.r.t. $P$. By Lemma A.8 and the fact that $E(\varphi|U_K) = 1$, there is a family of probability densities $(\varphi_n)_{n \in \mathbb{N}}$ and a family of indices in $K$, $(k_n)_{n \in \mathbb{N}}$, s.t. $\varphi_n$ is $\sigma\big(U_{k_n}\big)$-measurable and $\prod_{n \in \mathbb{N}_0} \varphi_n$ converges (unconditionally) in $L^1(P)$ and a.s. to $\varphi$. For $n \in \mathbb{N}$ define $(P_n, Q_n) = \Big((\prod_{m=1}^{n-1} \varphi) \cdot P, (\prod_{m=1}^{n} \varphi) \cdot P\Big)$. By the independence of $U$, we have $P_n, Q_n \in \triangle^\times$.

- We have $P_n(A|Z)(\omega) = Q_n(A|Z)(\omega)$ for a.e. $\omega \in \{K \overset{\text{a.s.}}{\subseteq} J\}$. Indeed, note that $P_n = \varphi_n \cdot Q_n$ and therefore $(P_n, Q_n) \in \triangle_{k_n}^{\times 2}$. By the definition of irrelevance, $P_n(A|Z)(\omega) = Q_n(A|Z)(\omega)$ for a.e. $\omega \in \{k_n \in J\} \overset{\text{a.s.}}{\supseteq} \{K \overset{\text{a.s.}}{\subseteq} J\}$
- Since $P_0 = P$, by induction on $n$, we have $\forall n \in \mathbb{N} : P(A|Z)(\omega) = Q_n(A|Z)(\omega)$ for a.e. $\omega \in \{K \overset{\text{a.s.}}{\subseteq} J\}$.
- Finally, $Q_n \to Q$ in $d_{L^1(P)}$ (c.f. Definition A.2). Let $C = \{K \overset{\text{a.s.}}{\subseteq} J\}$. By convergence $P(A, C|Z) \overset{\text{a.s.}}{=} Q_n(A, C|\ Z) \to Q(A, C|Z)$ in measure (Lemma A.9), and therefore $P(A, C|Z) \overset{\text{a.s.}}{=} Q(A, C|Z)$. □

The following corollary shows that irrelevance behaves well w.r.t. almost sure equality.

**Corollary 6.1.4**: Let $J$ be irrelevant to $X$ given $Z$ and let $K \overset{\text{a.s.}}{=} J$. Then $K$ is irrelevant to $X$ given $Z$.

*Proof*: Trivial. □

We want to show that an almost surely maximal irrelevant random index set exists. As we have seen in Lemma 6.2, the complement of the history is one irrelevant random index set. We will show in this section that it is actually maximal.

**Lemma 6.1.5**: An almost surely unique maximal irrelevance to $X$ given $Z$ exists.

*Proof*: Let $\mathfrak{I}$ denote the set of irrelevant random index sets for $X$ given $Z$. Set $M = \bigcup^{\text{a.s.}} \mathfrak{I}$. We show that $M$ is irrelevant to $X$ given $Z$. Then clearly, $M$ is a maximal irrelevant random index set. Let $i \in I$, $(P, Q) \in \triangle_i^{\times 2}$ and $A \in \sigma(X)$. We need to show that $P(A|Z)(\omega) = Q(A|Z)(\omega)$ for a.e. $\omega \in \{i \in M\}$. By the definition of almost sure union, there exists a sequence of irrelevant random index sets $(J_n)_{n \in \mathbb{N}}$ in $\mathfrak{I}$, s.t. $\{i \in M\} \overset{\text{a.s.}}{=} \bigcup_{n \in \mathbb{N}} \{i \in J_n\}$. Therefore it suffices to show $P(A|Z)(\omega) = Q(A|Z)(\omega)$. for a.e. $\omega \in \{i \in J_n\}$ and $n \in \mathbb{N}$. This is immediate by the definition of irrelevance and $J_n \in \mathfrak{I}$. $\square$

**Definition 6.1.6** (random index set of irrelevance): We call the almost surely unique maximal irrelevant random index set from Lemma 6.1.5 the *random index set of irrelevance to $X$ given $Z$*, or the *irrelevance to $X$ given $Z$* for short. We denote it by

$$\mathcal{I}(X|Z) := \overset{\text{a.s.}}{\bigcup} \{J : \Omega \to \mathfrak{P}(I) \mid J \text{ is irrelevant to } X \text{ given } Z\}.$$

Lemma 6.1.5 tells us that $\mathcal{I}(X|Z)$ is irrelevant to $X$ given $Z$. Note that $\mathcal{I}(X|Z)$ is only defined up to almost sure equality.

We will establish that $\mathcal{I}(X|Z) = \overline{\mathcal{H}(X|Z)}$, so it really is the dual to the history in this sense. Recall that our goal in this section is to prove the following direction of the fundamental theorem

$$\forall P \in \triangle^{\times} : X \perp\!\!\!\perp_P Y \mid Z \Rightarrow \mathcal{H}(X|Z) \cap \mathcal{H}(Y|Z) \overset{\text{a.s.}}{=} \emptyset.$$

Under the assumption that $\mathcal{I}(X|Z) = \overline{\mathcal{H}(X|Z)}$, this becomes equivalent to the dual notion

$$\forall P \in \triangle^{\times} : X \perp\!\!\!\perp_P Y \mid Z \Rightarrow \mathcal{I}(X|Z) \cup \mathcal{I}(Y|Z) \overset{\text{a.s.}}{=} I. \tag{3}$$

The reason why it is necessary to define $\mathcal{I}(X|Z)$ is that in the definition of $\mathcal{H}(X|Z)$ we have to show the measurability condition $\sigma(X) \in \sigma\left(U_{\mathcal{H}(X|Z)}, Z\right)$ and the independence $U_{\mathcal{H}(X|Z)} \perp\!\!\!\perp_P U_{\mathcal{H}(Y|Z)} \mid Z$ for all $P \in \triangle^{\times}$. This is not in direct connection with the provided assumption $\forall P \in \triangle^{\times} : X \perp\!\!\!\perp_P Y \mid Z$. The definition of $\mathcal{I}(X|Z)$, however, uses only probabilities in its definition. So there is a much more direct connection with this assumption that we can exploit.

We might get the idea to define the history directly through $\mathcal{H}(X|Z) := \overline{\mathcal{I}(X|Z)}$. But there are a few issues. Firstly, the definition of $\mathcal{I}(X|Z)$ is rather indirect and not very telling about what kind of $X$ and $Y$ fulfill $\mathcal{I}(X|Z) \cup \mathcal{I}(Y|Z) \overset{\text{a.s.}}{=} I$. Secondly, recall that it was relatively straightforward to prove the direction

$$\mathcal{H}(X|Z) \cap \mathcal{H}(Y|Z) \overset{\text{a.s.}}{=} \emptyset \Rightarrow \forall P \in \triangle^{\times} : X \perp\!\!\!\perp_P Y \mid Z \tag{4}$$

If we dualize with $\mathcal{I}(X|Z)$ this is equivalent to

$$\mathcal{I}(X|Z) \cup \mathcal{I}(Y|Z) \overset{\text{a.s.}}{=} I \Rightarrow \forall P \in \triangle^{\times} : X \perp\!\!\!\perp_P Y \mid Z \tag{5}$$

But now this statement is not easily proven directly. The way to prove (5) moves through showing $\mathcal{H}(X|Z) = \overline{\mathcal{I}(X|Z)}$ and dualizing back to (4).

In summary, $\mathcal{H}(X|Z)$ and $\mathcal{I}(X|Z)$ are dual notions that more suitable to prove the (4) and (3) directions of the fundamental theorem, respectively. The power of them is precisely that they are dual and come together to prove the fundamental theorem in full.

The next theorem proves that (3) holds in a weaker form. It states that for fixed $i \in I$ and $(P, Q) \in \triangle_i^{\times 2}$, that for a.e. $\omega$, morally, we either have that $i$ is irrelevant to $X$ given $Z$ or $i$ is irrelevant to $Y$ given $Z$ for this choice of $i$ and $(P, Q)$ and $\omega$.

**Theorem 6.1.7**: Let $\forall P \in \triangle^{\times} : X \perp\!\!\!\perp_P Y \mid Z$. Let $i \in I$, $(P, Q) \in \triangle_i^{\times 2}$ and $A \in \sigma(X), B \in \sigma(Y)$. Then

$$(P(A|Z) - Q(A|Z))(P(B|Z) - Q(B|Z)) \overset{\text{a.s.}}{=} 0.$$

Therefore, for a.e. $\omega \in \Omega$, $P(A|Z)(\omega) = Q(A|Z)(\omega)$ or $P(B|Z)(\omega) = Q(A|Z)(\omega)$.

*Proof*: Let $f : \Omega \to \mathbb{R}_{>0}$ $\sigma(U_i)$-measurable s.t. $P = f \cdot Q$. We use the fact that for the positive $\sigma(U_i)$-measurable density $g = \frac{1+f}{2}$, we can define $R := \frac{P+Q}{2} = g \cdot Q \in \triangle^{\times}$. Therefore $X \perp\!\!\!\perp_R Y \mid Z$ by assumption. Denoting $E$ the expectation w.r.t $Q$ and noting that $g > 0$,

$$R(A|Z)R(B|Z) \overset{\text{a.s.}}{=} R(A,B|Z)$$

$$\Leftrightarrow E(g1_A|Z)E(g1_B|Z) \overset{\text{a.s.}}{=} E(g|Z)E(g1_A 1_B|Z)$$

We can now substitute $g$ by $\frac{f+1}{2}$ and multiply by 4.

$$\Leftrightarrow (E(f1_A|Z) + E(1_A|Z))(E(f1_B|Z)E(1_B|Z)) \overset{\text{a.s.}}{=} (E(f|Z) + 1)(E(f1_A 1_B|Z) + E(1_A 1_B|Z))$$

By the independence of $X$ and $Y$ given $Z$ w.r.t. $P$ and $Q$, we have $E(f1_A|Z)E(f1_B|Z) \overset{\text{a.s.}}{=} E(f|Z)E(f1_A 1_B|Z)$ and $E(1_A|Z)E(1_B|Z) \overset{\text{a.s.}}{=} E(1_A 1_B|Z)$ respectively. Multiplying out and canceling these terms, we get

$$\Leftrightarrow E(f1_A|Z)E(1_B|Z) + E(1_A|Z)E(f1_B|Z) \overset{\text{a.s.}}{=} E(f|Z)E(1_A 1_B|Z) + E(f1_A 1_B|Z)$$

Dividing by $E(f|Z)$, we get

$$\Leftrightarrow P(A|Z)Q(B|Z) + Q(A|Z)P(B|Z) \overset{\text{a.s.}}{=} Q(A,B|Z) + P(A,B|Z)$$

Using the independence w.r.t. $P$ and $Q$ again,

$$\Leftrightarrow P(A|Z)Q(B|Z) + Q(A|Z)P(B|Z) \overset{\text{a.s.}}{=} Q(A|Z)Q(B|Z) + P(A|Z)P(B|Z)$$

which we can factorize as

$$\Leftrightarrow (P(A|Z) - Q(A|Z))(P(B|Z) - Q(B|Z)) \overset{\text{a.s.}}{=} 0.$$

□

## 6.2 The fundamental theorem for the random index set of irrelevance

**Theorem 6.2.1**: Let $i \in I$. If the set

$$\triangle_{i,C}^{\times 2} := \left\{(P,Q) \in \triangle_i^{\times 2} : P(A|Z)(\omega) \neq Q(A|Z)(\omega) \text{ for a.e. } \omega \in C\right\}$$

is not empty, then it is dense w.r.t. the product metric

$$(d_1 \times d_1)((P,Q),(P',Q')) := d_1(P,P') + d_1(Q,Q') \text{ (c.f. Definition A.2)}$$

in $\triangle_i^{\times 2} \subseteq (\triangle^{\times})^2$.

*Proof*: Let $(P,Q) \in \triangle_{i,C}^{\times 2}$, i.e. $P(A|Z)(\omega) \neq Q(A|Z)(\omega)$ for a.e. $\omega$ in $C$. Since $(P,Q) \in \triangle_i^{\times 2}$, $Q = q \cdot P$ for a $\sigma(U_i)$-measurable probability density $q$.

Now let $(P',Q') \in \triangle_i^{\times 2}$. We need to show that $(P',Q')$ is in the $(d_1 \times d_1)$-closure of $\triangle_{i,C}^{\times 2}$.

By Lemma A.8, there is a family of densities $(\varphi_n)_{n \in \mathbb{N}_0}$ and indices $\{i_n\}_{n \in \mathbb{N}}$, s.t.

- $\mathbb{E}(\varphi_0|U) \overset{\text{a.s.}}{=} 1$.
- $\forall n \in \mathbb{N} : \varphi_n$ is $\sigma\left(U_{i_n}\right)$-measurable.
- $\prod_{n \in \mathbb{N}_0} \varphi_n$ converges (unconditionally) in $L^1(P)$ and a.s. pointwise to $\frac{\mathrm{d}P'}{\mathrm{d}P}$.

W.l.o.g. we can assume $i = i_1$. Then there is a $\sigma(U_i)$-measurable probability density $q'$, s.t. $Q' = q' \cdot P'$. Set $\varphi_1' = q' \cdot \varphi_1$ and for $n \in \mathbb{N}_0 \setminus \{1\}$, set $\varphi_n' = \varphi_n$. Then $\prod_{n \in \mathbb{N}_0} \varphi_n' = \frac{\mathrm{d}Q'}{\mathrm{d}P}$, since for any $A \in \mathcal{A}$, $\int_A \frac{\mathrm{d}Q'}{\mathrm{d}P} \,\mathrm{d}P = \int_A \mathrm{d}Q' = \int_A q' \,\mathrm{d}P' = \int_A q' \prod_{n \in \mathbb{N}_0} \varphi_n \,\mathrm{d}P = \int_A \prod_{n \in \mathbb{N}_0} \varphi_n' \,\mathrm{d}P$.

We now show that $(P',Q')$ is in the $(d_1 \times d_1)$-closure of $\triangle^{\times}(i,C)$. The idea is to connect $(P,Q)$ to $(P',Q')$ in a well-behaved curve $\lambda \mapsto (R_\lambda, R_\lambda')$. The condition that a point on this curve is in $\triangle^{\times}(i,C)$ corresponds to a family of polynomials evaluated at $\lambda$ being almost surely unequal to zero. Since the starting point $(P,Q) = (R_0, R_0')$ is in $\triangle_{i,C}^{\times 2}$, almost all the polynomials cannot be the zero polynomial. Therefore the polynomials are almost surely unequal to zero at almost all $\lambda$. Therefore almost all points on the curve are in $\triangle_{i,C}^{\times 2}$.

1. Finite dimensional case[2]: Suppose $\exists m \in \mathbb{N} : \forall n \geq m : \varphi_n = 1$. Then define

[2]I acknowledge Scott Garrabrant, who simplified the proof of the finite dimensional case, in the finite case [5], from an induction to one step. These ideas are used to simplify the proof here.

$$
\begin{aligned}
p : [0,1] \times \mathbb{R}^{\{0,..,m\}} &\to [0,1] \\
(\lambda, x) &\mapsto \prod_{n=0}^{m} \big(\lambda x_n + \overline{\lambda}\big).
\end{aligned}
$$

Note that $p(\cdot, x)$ is a polynomial for $x \in \mathbb{R}^{\{0,\dots,m\}}$. Set $\Phi = (\varphi_n)_{n=1}^{m}$ and $\Phi' = (\varphi_n')_{n=1}^{m}$. Define $\varphi_\lambda(\omega) = p(\lambda, \Phi(\omega))$ and $\varphi'_\lambda(\omega) = p(\lambda, \Phi'(\omega))$. The following properties hold for $\varphi_\lambda$ and $\varphi'_\lambda$.

(i) We claim that $\varphi_\lambda$ and $\varphi'_\lambda$ are probability densities w.r.t. $P$. Indeed,

$$
\begin{aligned}
\int \varphi_\lambda \,\mathrm{d}P &= \int \prod_{n=0}^{m} \big(\lambda \varphi_n + \overline{\lambda}\big) \,\mathrm{d}P \\
&= \int \mathbb{E}\big(\lambda \varphi_0 + \overline{\lambda} | U\big) \prod_{n=1}^{m} \big((\lambda \varphi_n) + \overline{\lambda}\big) \,\mathrm{d}P \\
&= \prod_{n=1}^{m} \int \big((\lambda \varphi_n) + \overline{\lambda}\big) \,\mathrm{d}P = 1,
\end{aligned}
$$

because of $\mathbb{E}(\varphi_0 | U) \overset{\text{a.s.}}{=} 1$ and the independence of $U$. $\varphi'_\lambda$ is proven analogously.

(ii) We claim that $\lambda \mapsto \varphi_\lambda$ and $\lambda \mapsto \varphi'_\lambda$ are continuous maps into $L^1(P)$. Indeed, let $\lambda, \lambda' \in [0,1]$. Then, since $\lambda \varphi_n - \overline{\lambda} - \big(\lambda' \varphi_n - \overline{\lambda'}\big) = (\lambda - \lambda')(\varphi_n - 1)$,

$$
\begin{aligned}
\|\varphi_\lambda - \varphi_{\lambda'}\|_{L^1(P)} &= \int \left| \prod_{n=0}^{m} \big(\lambda \varphi_n - \overline{\lambda}\big) - \prod_{n=0}^{m} \big(\lambda' \varphi_n - \overline{\lambda'}\big) \right| \mathrm{d}P \\
&\le \|\mathbb{E}(|(\lambda - \lambda')(\varphi_0 - 1)| \; |U)\|_\infty \prod_{n=1}^{m} \|(\lambda - \lambda')(\varphi_n - 1)\|_{L^1(P)} \\
&= (2(\lambda - \lambda'))^{m+1}.
\end{aligned}
$$

Set $R_\lambda := \varphi_\lambda \cdot P$ and $R'_\lambda := \varphi'_\lambda \cdot P$. Now $\lambda \in [0,1] \mapsto (R_\lambda, R'_\lambda)$ is a continuous curve connecting $(P, Q)$ and $(P', Q')$. Our goal is to show that for $\lambda$ near 1, we have $R_\lambda(A|Z) \neq R'_\lambda(B|Z)$ a.s. on $C$. Set $E_{\Phi,A} = (\mathbb{E}(\varphi_n 1_A | Z))_{n=1}^{m}$, $E_{\Phi',A} = (\mathbb{E}(\varphi'_n 1_A | Z))_{n=1}^{m}$. $E_\Phi = (\mathbb{E}(\varphi_n | Z))_{n=1}^{m}$, $E_{\Phi'} = (\mathbb{E}(\varphi'_n | Z))_{n=1}^{m}$.

Now note that for a.e. $\omega \in \Omega$ conditional expectation is linear, and so

$$
\begin{aligned}
& R_\lambda(A|Z)(\omega) = R'_\lambda(A|Z)(\omega) \\
\Leftrightarrow\; & \mathbb{E}(\varphi_\lambda 1_A | Z) \mathbb{E}(\varphi'_\lambda | Z)(\omega) = \mathbb{E}(\varphi'_\lambda 1_A | Z)(\omega) \mathbb{E}(\varphi_\lambda | Z)(\omega) \\
\Leftrightarrow\; & \mathbb{E}(p(\lambda, \Phi) 1_A | Z) \mathbb{E}(p(\lambda, \Phi') | Z)(\omega) = \mathbb{E}(p(\lambda, \Phi') 1_A | Z)(\omega) \mathbb{E}(p(\lambda, \Phi) | Z)(\omega) \\
\Leftrightarrow\; & \underbrace{p\big(\lambda, E_{\Phi,A}(\omega)\big) p(\lambda, E_{\Phi'}(\omega)) - p\big(\lambda, E_{\Phi',A}(\omega)\big) p(\lambda, E_\Phi(\omega))}_{=: p'(\lambda, \omega)} = 0
\end{aligned}
$$

Let $p'(\lambda, \omega) = p\big(\lambda, E_{\Phi,A}(\omega)\big) p(\lambda, E_{\Phi'}(\omega)) - p\big(\lambda, E_{\Phi',A}(\omega)\big) p(\lambda, E_\Phi(\omega))$. Clearly, $p' : [0,1] \times \Omega \to \mathbb{R}$ is measurable and $p'(\cdot, \omega)$ is a polynomial for all $\omega \in \Omega$. Then we have just proved that $R_\lambda(A|Z)(\omega) = R'_\lambda(A|Z)(\omega) \Leftrightarrow p'(\lambda, \omega) = 0$. Therefore, because $(R_0, R'_0) = (P, Q) \in \triangle^{\times 2}_{i,C}$, we have $p'(0, \omega) \neq 0$ for a.e. $\omega \in C$. Therefore $p'(\cdot, \omega)$ is not the zero polynomial for a.e. $\omega \in C$ and the set $\{\lambda \in [0,1] : p'(\lambda, \omega) = 0\}$ is finite for a.e. $\omega \in C$ and therefore a nullset w.r.t. to the Lebesgue measure $\mathcal{L}$ on $[0,1]$. By Fubini, $(\mathcal{L} \times P)\{(\lambda, \omega) \in [0,1] \times C : p'(\lambda, \omega) = 0\} = 0$. Again, by Fubini, for $\mathcal{L}$-a.e. $\lambda \in [0,1]$, we have $P(\{\omega \in C : p'(\lambda, \omega) = 0\}) = 0$. So we can choose a sequence $\lambda_n \in [0,1]$, s.t. $\lambda_n \to 1$ and for $n \in \mathbb{N}$, $p'(\lambda_n, \omega) \neq 0$ for a.e. $\omega \in C$. Then by construction, $\big(R_{\lambda_n}, R'_{\lambda_n}\big) \in \triangle^{\times 2}_{i,C}$, while $R_{\lambda_n} \to P'$ and $R'_{\lambda_n} \to Q'$ in $d_1$. Therefore, $(P', Q')$ is in the $(d_1 \times d_1)$-closure of $\triangle^{\times 2}_{i,C}$.

2. General case: The key idea is that we can truncate the infinite product high enough to get a finite product that is arbitrarily close to the infinite product and then use the finite dimensional case.

   Let $\varepsilon > 0$. Choose, $m \in \mathbb{N}$ s.t.
   $$\left\| \prod_{n=0}^{m} \varphi_n - \varphi \right\|_{L^1(P)} < \varepsilon \quad \text{and} \quad \left\| \prod_{n=0}^{m} \varphi'_n - \varphi' \right\|_{L^1(P)} < \varepsilon.$$
   Set $\varphi_{\le m} = \prod_{n=0}^{m} \varphi_n$ and $\varphi'_{\le m} = \prod_{n=0}^{m} \varphi'_n$. Set $P'_{\le m} = \varphi_{\le m} \cdot P$ and $Q'_{\le m} = \varphi'_{\le m} \cdot P$. Then by the choice of $m$, $d_1(P'_{\le m}, P') < \varepsilon$ and $d_1(Q'_{\le m}, Q') < \varepsilon$. Clearly, $(P'_{\le m}, Q'_{\le m})$ fulfills the assumption of finite dimensional case and there is $(P'', Q'') \in \triangle^{\times 2}_{i,C}$ s.t. $(d_1 \times d_1)((P'', Q''), (P'_{\le m}, Q'_{\le m})) < \varepsilon$. By the triangle inequality, $(d_1 \times d_1)((P'', Q''), (P', Q')) < 3\varepsilon$. Since $\varepsilon > 0$ was arbitrary, we have shown that $(P', Q')$ is in the $(d_1 \times d_1)$-closure of $\triangle^{\times 2}_{i,C}$. □

We can now use the theorems above to show one direction of the fundamental theorem for the random index set of irrelevance.

**Theorem 6.2.2**: If $\forall P \in \triangle^{\times} : X \perp\!\!\!\perp_P Y \mid Z$, then $\mathcal{J}(X|Z) \cup \mathcal{J}(Y|Z) \overset{\text{a.s.}}{=} I$.

*Proof*: Let $i \in I$, we show that $\{i \notin \mathcal{J}(X|Z)\} \overset{\text{a.s.}}{\subseteq} \{i \in \mathcal{J}(Y|Z)\}$. Suppose not, then $C_0 := \{i \notin \mathcal{J}(X|Z)\} \setminus \{i \in \mathcal{J}(Y|Z)\} \notin \mathcal{N}$. Then by the definition of irrelevance, there is $A \in \sigma(X)$ and $(P, Q) \in \triangle^{\times 2}_{i}$ and a $C \in \sigma(Z) \setminus \mathcal{N}$ s.t. $C \overset{\text{a.s.}}{\subseteq} C_0$ and $P(A|Z)(\omega) \neq Q(A|Z)(\omega)$ for a.e. $\omega \in C$.

Towards a contradiction, we show that $C \overset{\text{a.s.}}{\subseteq} \{i \in \mathcal{J}(Y|Z)\}$. For this, let $(P', Q') \in \triangle^{\times 2}_{i}$ and $B \in \sigma(Y)$. It suffices to show that $P'(B|Z)(\omega) = Q'(B|Z)(\omega)$ for a.e. $\omega \in C$, since then the random index set $J(\omega) := \{i\}$ if $\omega \in C$ and $\emptyset$ else; is irrelevant to $Y$ given $Z$. For this, it suffices to show $P'(B, C|Z) \overset{\text{a.s.}}{=} Q'(B, C|Z)$, since $C \in \sigma(Z)$.

We use Theorem 6.2.1. Since $(P, Q) \in \triangle^{\times 2}_{i,C}$, $\triangle^{\times 2}_{i,C}$ is not empty and therefore dense in $\triangle^{\times}$ $(i)$ w.r.t. $d_1 \times d_1$. Let $(P_n, Q_n) \in \triangle^{\times 2}_{i,C}$ s.t. $(P_n, Q_n) \to (P', Q')$. By definition of $\triangle^{\times 2}_{i,C}$, we have for all $n \in \mathbb{N}$ that $P_n(A|Z)(\omega) \neq Q_n(A|Z)(\omega)$ for a.e. $\omega \in C$ and therefore by Theorem 6.1.7, $P_n(B|Z)(\omega) = Q_n(B|Z)(\omega)$ for a.e. $\omega \in C$. Then by Lemma A.9, $P_n(B, C|Z) \to P'(B, C|Z)$ in measure, while $P_n(B, C|Z) \overset{\text{a.s.}}{=} Q_n(B, C|Z) \to Q'(B, C|Z)$. Since a limit in measure is almost surely unique, $P'(B, C|Z) \overset{\text{a.s.}}{=} Q'(B, C|Z)$. □

## 6.3 The duality between history and irrelevance

The goal of this subsection is to prove that $\mathcal{H}(X|Z) = \overline{\mathcal{J}(X|Z)}$. We first prove that $\mathcal{J}(X|Z)$ behaves like $\mathcal{H}(X|Z)$ w.r.t. independence. Firstly, it is easy to see that $\mathcal{H}\left(U_{\overline{\mathcal{H}(X|Z)}}\right) = \overline{\mathcal{H}(X|Z)}$, so $X \perp U_{\overline{\mathcal{H}(X|Z)}}$. This translates to $X \perp U_{\mathcal{J}(X|Z)}$. Morally, this means that everything irrelevant to $X$ given $Z$ is independent of $X$ given $Z$.

**Theorem 6.3.1**: We have $\forall P \in \triangle^{\times} : X \perp\!\!\!\perp_P U_{\mathcal{J}(X|Z)} \mid Z$.

*Proof*: Let $P \in \triangle^{\times}$ and $A \in \sigma(X, Z)$. Let $K \subseteq I$ finite and $B = \bigtimes_{k \in K} B_k$, where $B_k \subseteq \mathrm{Val}(U_k)$ measurable. By Corollary 4.6, it suffices to show that $A \perp\!\!\!\perp_P \{K \subseteq \mathcal{J}(X|Z), U_K = B\} \mid Z$, since these sets form a $\cap$-stable generator of $\sigma\left(U_{\mathcal{J}(X|Z)}\right)$. Since $C := \{K \subseteq \mathcal{J}(X|Z)\} \in \sigma(Z)$, it suffices to show that $A \cap C \perp\!\!\!\perp_P U_K \mid Z$. Therefore, w.l.o.g. we can assume $A \subseteq C$. Let $\varphi_k$ be arbitrary, positive, $\sigma(U_k)$-measurable probability densities w.r.t. $P$. Set $\varphi = \prod_{k \in K} \varphi_k$ and $Q = \varphi \cdot P$. Then $(P, Q) \in \triangle^{\times 2}_{K}$ and by Lemma 6.1.3, $P(A|Z) \overset{\text{a.s.}}{=} Q(A|Z)$ because $A \subseteq C = \{K \subseteq \mathcal{J}(X|Z)\}$.

We now show that $P(A|Z) \overset{\text{a.s.}}{=} P(A|U_K, Z)$, which is equivalent to $A \perp\!\!\!\perp_P U_K \mid Z$. Now let $D \in \sigma(Z)$. Let $E$ denote the expectation w.r.t. $P$.

$$
\begin{aligned}
\int_D P(A|Z)\,\mathrm{d}Q &= \int_D Q(A|Z)\,\mathrm{d}Q \\
&= \int_D 1_A\,\mathrm{d}Q \\
&= \int_D f 1_A\,\mathrm{d}P \\
&= \int_D E(f 1_A|U_K, Z)\,\mathrm{d}P \\
&= \int_D P(A|U_K, Z)\,\mathrm{d}Q
\end{aligned}
$$

Therefore, taking conditional expectation w.r.t. $U_K$,

$$
\int \varphi E(1_D(P(A|Z) - P(A|U_K, Z))|U_K)\,\mathrm{d}P = 0
$$

Clearly, we can approximate the indicator variable of any rectangle $B = \bigcap_{k \in K} B_k$, where $B_k \in \sigma(U_k)$, uniformly by the positive probability density $\varphi^\varepsilon = \prod_{k \in K} \varphi_k^\varepsilon$, where $\varphi_k^\varepsilon = \frac{1_B + \varepsilon}{P(B) + \varepsilon}$. As $\varepsilon \to 0$, $\varphi^\varepsilon \to 1_B$ uniformly. Therefore,

$$
\int_B E(1_D(P(A|Z) - P(A|U_K, Z))|U_K)\,\mathrm{d}P = 0 \tag{6}
$$

for all these choices of $B$. This set of rectangles is a $\cap$-stable generator of $\sigma(U_K)$. The set of $B \subseteq \mathrm{Val}(X_K)$ which fulfill (6), clearly form a Dynkin system, therefore (6) holds for all $B \in \sigma(U_K)$ measurable.

Therefore, since $E(1_D(P(A|Z) - P(A|U_K, Z))|U_K)$ is $\sigma(U_K)$-measurable, we have

$$
E(1_D(P(A|Z) - P(A|U_K, Z))|U_K) \overset{\text{a.s.}}{=} 0.
$$

Now let $F \in \sigma(U_K)$ then $E(1_{D \cap F}(P(A|Z) - P(A|U_K, Z))|U_K) \overset{\text{a.s.}}{=} 0$. Then

$$
\int_{D \cap F} P(A|Z) - P(A|U_K, Z)\,\mathrm{d}P = 0 \tag{7}
$$

for all $D \in \sigma(Z)$ and $F \in \sigma(U_K)$. Since $\{D \cap F : D \in \sigma(Z)\}$ is a $\cap$-stable generator of $\sigma(U_K, Z)$ and the set of all $G \in \sigma(U_K, Z)$ for which (7) holds, form a Dynkin system, (7) holds for all $G \in \sigma(U_K, Z)$. Since $P(A|Z) - P(A|U_K, Z)$ is $\sigma(U_K, Z)$-measurable, we have $P(A|Z) \overset{\text{a.s.}}{=} P(A|U_K, Z)$. □

Next, we use this result to see that $\mathcal{I}(X|Z)$ disintegrates $X$ given $Z$, just like the history $\mathcal{H}(X|Z)$.

**Lemma 6.3.2**: We have $\forall P \in \triangle^\times : U_{\overline{\mathcal{I}(X|Z)}} \perp\!\!\!\perp_P U_{\mathcal{I}(X|Z)} \mid Z$.

*Proof*: By Theorem 6.3.1, $\forall P \in \triangle^\times : X \perp\!\!\!\perp_P U_{\mathcal{I}(X|Z)} \mid Z$. Let $B \in \sigma\big(U_{\mathcal{I}(X|Z)}\big)$. Then $\forall P \in \triangle^\times$ $X \perp\!\!\!\perp_P B \mid Z$, and therefore by Theorem 6.2.2, $\mathcal{I}(X|Z) \cup \mathcal{I}(B|Z) \overset{\text{a.s.}}{=} I$ and thus $\mathcal{I}(B|Z) \supseteq \overline{\mathcal{I}(X|Z)}$. By Theorem 6.3.1, $\forall P \in \triangle^\times : B \perp\!\!\!\perp_P U_{\mathcal{I}(B|Z)} \mid Z$ and by Lemma 4.3.6, $\sigma\big(U_{\mathcal{I}(B|Z)}, Z\big) \supseteq \sigma\big(U_{\overline{\mathcal{I}(X|Z)}}, Z\big)$. Therefore, $\forall P \in \triangle^\times : B \perp\!\!\!\perp_P U_{\overline{\mathcal{I}(X|Z)}}$. Since $B \in \sigma\big(U_{\mathcal{I}(X|Z)}\big)$ was arbitrary, we have $\forall P \in \triangle^\times : U_{\overline{\mathcal{I}(X|Z)}} \perp\!\!\!\perp U_{\mathcal{I}(X|Z)} \mid Z$. □

We now show that the independence in the last two lemmas is compositional, so that we can combine the two independence statements into one.

**Theorem 6.3.3**: We have $\forall P \in \triangle^\times : \left(X, U_{\overline{\mathcal{J}(X|Z)}}\right) \perp\!\!\!\perp_P U_{\mathcal{J}(X|Z)} \mid Z$.

*Proof*: Let $P \in \triangle^\times$. Let $K \subseteq I$ be finite and nonempty. Let $C = \left\{K \subseteq \overline{\mathcal{J}(X|Z)}\right\}$. By Corollary 4.6 it suffices to show $\sigma(X, Z, U_K)|_C \perp\!\!\!\perp_P B \mid Z$ for all $B \in \sigma\left(U_{\mathcal{J}(X|Z)}, Z\right)$.

Let $A \in \sigma(X, Z)|_C$, let $B \in \sigma\left(U_{\mathcal{J}(X|Z)}, Z\right)$. Let $D = \bigcap_{k \in K} D_k$, where $D_k \in \sigma(U_k)$. It suffices to show that $A \cap D \perp\!\!\!\perp_P B \mid Z$, since sets of the form $A \cap D$ form a $\cap$-stable generator of $\sigma(X, Z, U_K)|_C$. W.l.o.g. we can assume $B \subseteq C$, since $A \subseteq C \in \sigma(Z)$.

We will approximate $1_D$ by densities, so let $\varphi$ be an arbitrary positive $\sigma(U_K)$-measurable probability density w.r.t. $P$ and set $Q = \varphi \cdot P$. Then $Q \in \triangle^\times$ and by Theorem 6.3.1 $X \perp\!\!\!\perp_Q U_{\mathcal{J}(X|Z)} \mid Z$ and therefore $Q(A|Z)Q(B|Z) \overset{\text{a.s.}}{=} Q(A, B|Z)$. Let $E$ be the expectation w.r.t $P$. Then

$$E(\varphi 1_A|Z)E(\varphi 1_B|Z) \overset{\text{a.s.}}{=} E(\varphi|Z)E(\varphi 1_A 1_B|Z). \tag{8}$$

Since $1_C\varphi$ is $\sigma(U_K, Z)|_C \subseteq \sigma\left(U_{\overline{\mathcal{J}(X|Z)}}, Z\right)|_C$-measurable, we have $1_C\varphi \perp\!\!\!\perp_Q B \mid Z$ by Lemma 6.3.2. Since $B \subseteq C$, we have $\varphi \perp\!\!\!\perp_Q B \mid Z$. Applying this in (8), we have

$$\begin{aligned} E(\varphi 1_A|Z)E(\varphi|Z)E(1_B|Z) &\overset{\text{a.s.}}{=} E(\varphi|Z)E(\varphi 1_A 1_B|Z) \\ \Leftrightarrow E(\varphi 1_A|Z)E(1_B|Z) &\overset{\text{a.s.}}{=} E(\varphi 1_A 1_B|Z) \end{aligned}$$

By choosing $\varphi^\varepsilon = \prod_{n=1}^k \varphi_k^\varepsilon$, where $\varphi_\varepsilon^k = \frac{D_k+\varepsilon}{P(D_k)+\varepsilon}$, we have $\varphi^\varepsilon \to 1_D$ uniformly and therefore

$$\begin{aligned} E(1_D 1_A|Z)E(1_B|Z) &\overset{\text{a.s.}}{=} E(1_D 1_A 1_B|Z) \\ \Leftrightarrow A \cap D &\perp\!\!\!\perp_P B \mid Z \end{aligned}$$

In conclusion, we have $\left(X, U_{\overline{\mathcal{J}(X|Z)}}\right) \perp\!\!\!\perp_P U_{\mathcal{J}(X|Z)} \mid Z$. □

**Theorem 6.3.4**: $\mathcal{H}(X|Z) \overset{\text{a.s.}}{=} \overline{\mathcal{J}(X|Z)}$.

*Proof*: We show both inclusions.
'$\overset{\text{a.s.}}{\supseteq}$': It suffices to show $\overline{\mathcal{H}(X|Z)} \subseteq \mathcal{J}(X|Z)$. This follows immediately from showing that $\overline{\mathcal{H}(X|Z)}$ is irrelevant to $X$ given $Z$. We have shown this in Lemma 6.2.
'$\overset{\text{a.s.}}{\subseteq}$': It suffices to show that $\overline{\mathcal{J}(X|Z)}$ generates $X$ given $Z$. $\overline{\mathcal{J}(X|Z)}$ disintegrates $X$ given $Z$ by Lemma 6.3.2. Let $A \in \sigma(X)$ and $H = \overline{\mathcal{J}(X|Z)}$. It is left to show that $A \in \sigma(U_H, Z)$. Let $P \in \triangle^\times$ with expectation $E$. It suffices to show that $1_A \overset{\text{a.s.}}{=} E(1_A|U_H, Z)$.

By Theorem 6.3.3, we have $(X, U_H) \perp\!\!\!\perp_P U_{\overline{H}} \mid Z$. We first show that for all $C \in \sigma(U)$, we have

$$E(1_C 1_A|Z) \overset{\text{a.s.}}{=} E(1_C E(1_A|U_H, Z)|\ Z) \tag{9}$$

For this, it suffice to show (9) for all $C = B \cap D$, where $B \in \sigma(U_{\overline{H}})$ and $D \in \sigma(U_H)$, since sets of this form are a $\cap$-stable generator of $\sigma(U_{\overline{H}}, U_H) = \sigma(U)$ (Lemma 4.3.6) and the sets $C$ that fulfill (9) clearly form a Dynkin system. Now

$$\begin{aligned} & E(1_B 1_D 1_A|Z) \overset{\text{a.s.}}{=} E(1_B 1_D E(1_A|U_H, Z)|Z) \\ \Leftrightarrow\ & E(1_B|Z)E(1_D 1_A|Z) \overset{\text{a.s.}}{=} E(1_B|Z)E(E(1_D 1_A|U_H, Z)|Z) \\ \Leftarrow\ & E(1_D 1_A|Z) \overset{\text{a.s.}}{=} E(E(1_D 1_A|U_H, Z)|Z) \\ \Leftrightarrow\ & E(1_D 1_A|Z) \overset{\text{a.s.}}{=} E(1_D 1_A|Z) \end{aligned}$$

Finally, integrating (9), we have $\int_C 1_A \,\mathrm{d}P = \int_C E(1_A|U_H, Z) \,\mathrm{d}P$ for all $C \in \sigma(U)$. Since we have generally assumed that $X$ and $Z$ are $\sigma(U)$-measurable, it follows that $1_A \overset{\text{a.s.}}{=} E(1_A|U_H, Z)$. □

**Theorem 6.3.5** (completeness of structural independence): Structural independence is complete: If independence holds in all product probability distributions $P \in \triangle^\times$, then the independence is structural. More precisely, $\forall P \in \triangle^\times : X \perp\!\!\!\perp_P Y \mid Z \Rightarrow X \perp Y \mid Z$.

*Proof*: This follows from Theorem 6.2.2 and Theorem 6.3.4. □

## 6.4 The fundamental theorem

We now state the fundamental theorem of structural independence with all its assumptions.

**Theorem 6.4.1** (the fundamental theorem of structural independence): Let $U = (U_i)_{i\in I}$ be an independent family of random elements on a probability space $(\Omega, \mathcal{A}, \mathbb{P})$. Let $X, Y$ and $Z$ be $\sigma(U)$-measurable random elements.

Morally, the fundamental theorem states that, in general, we can conclude that $X \perp\!\!\!\perp_{\mathbb{P}} Y \mid Z$ using only the fact that $U$ is independent w.r.t. $\mathbb{P}$ if and only if $\mathcal{H}(X|Z) \cap \mathcal{H}(Y|Z) = \emptyset$ $\mathbb{P}$-a.s.

More formally, let

$$\begin{aligned}\triangle^{\times} = \{P : \mathcal{A} \to \mathbb{R} \mid\ & P \text{ is a probability distributions} \\ & \text{and } U \text{ is independent w.r.t. } P \\ & \text{and } P \sim \mathbb{P}\}.\end{aligned}$$

Then $\forall P \in \triangle^{\times} : X \perp\!\!\!\perp_P Y \mid Z \Leftrightarrow \mathcal{H}(X|Z) \cap \mathcal{H}(Y|Z) = \emptyset$ $\mathbb{P}$-a.s. With Definition 5.12, this can also be written as

$$\forall P \in \triangle^{\times} : X \perp\!\!\!\perp_P Y \mid Z \Leftrightarrow X \perp Y \mid Z$$

*Proof*: '$\Rightarrow$' is Theorem 6.3.5 '$\Leftarrow$' is Theorem 6.1. □

Furthermore, we can show that if structural independence does not hold, then independence is 'rare'. The setting is the same as in Theorem 6.4.1. Recall the definition of the metric $d_1$ on $\triangle^{\times}$ from Definition A.2.

**Lemma 6.4.2**: $\{P \in \triangle^{\times} : X \perp\!\!\!\perp_P Y \mid Z\}$ is closed in the topology induced by $d_1$.

*Proof*: Let $P_n$ be a sequence in $\{P : X \perp\!\!\!\perp_P Y \mid Z\}$ s.t. $P_n \xrightarrow{d_1} P$. By Lemma A.9, $P_n(A|Z) \to P(A|Z)$ in $P$-measure for all $A \in \mathcal{A}$. Let $A \in \sigma(X)$ and $B \in \sigma(Y)$. Then

$$P(A|Z)P(B|Z) \leftarrow P_n(A|Z)P_n(B|Z) \overset{\text{a.s.}}{=} P_n(A, B|Z) \to P(A, B|Z).$$

□

The following theorem is similar in spirit to Theorem 6.2.1.

**Theorem 6.4.3**: Suppose there is $\mathcal{D} \subseteq \triangle^{\times}$ open in the topology induced by $d_1$, s.t. $\forall P \in \mathcal{D} : X \perp\!\!\!\perp_P Y \mid Z$. Then $X \perp Y \mid Z$.

*Proof*: Let $P \in \mathcal{D}$ and $\varepsilon > 0$ s.t. $B_\varepsilon(P) \coloneqq \{Q \in \triangle^{\times} : d_1(P, Q) < \varepsilon\} \subseteq \mathcal{D}$

Let $Q \in \triangle^{\times}$. We need to show $X \perp\!\!\!\perp_Q Y \mid Z$.

Since $P \sim Q$, by Lemma A.8, there is a family of densities $(\varphi_n)_{n\in\mathbb{N}_0}$ and indices $\{i_n\}_{n\in\mathbb{N}}$, s.t.

- $\mathbb{E}(\varphi_0|U) \overset{\text{a.s.}}{=} 1$.
- $\forall n \in \mathbb{N} : \varphi_n$ is $\sigma\left(U_{i_n}\right)$-measurable.
- $\prod_{n\in\mathbb{N}_0} \varphi_n$ converges (unconditionally) in $L^1(P)$ and a.s. pointwise to $\frac{\mathrm{d}Q}{\mathrm{d}P}$.

Set $Q_{\leq n} \coloneqq \prod_{k=0}^{n} \varphi_k \cdot Q$. Then $Q_{\leq n} \xrightarrow{d_1} Q$. Therefore it suffices to show that $\forall n \in \mathbb{N} : X \perp\!\!\!\perp_{\{Q_{\leq n}\}} Y \mid Z$.

Let $n \in \mathbb{N}$. Define $p : [0, 1] \times \Omega \to \mathbb{R}; (\lambda, \omega) \mapsto \prod_{k=0}^{n} \left(\lambda\varphi_k + \overline{\lambda}\right)$. Set $\Phi = (\varphi_k)_{k=1}^{n}$ and set $\varphi_\lambda(\omega) = p(\lambda, \omega)$. Set $R_\lambda \coloneqq \varphi_\lambda \cdot P$. Like in the proof of Theorem 6.2.1, we can see that $\varphi_\lambda$ is a probability density w.r.t. $P$ for all $\lambda \in [0, 1]$. Furthermore, the map $\lambda \mapsto R_\lambda$ is continuous in $d_1$.

Let $A \in \sigma(X)$ and $B \in \sigma(Y)$. For $\lambda \in [0, 1]$, we have

$$A \perp\!\!\!\perp_{\{R_\lambda\}} B \mid Z \Leftrightarrow R_\lambda(A|Z)R_\lambda(B|Z) \overset{\text{a.s.}}{=} R_\lambda(A,B|Z)$$

$$\Leftrightarrow \underbrace{\left(E(\varphi_\lambda 1_A|Z) + E\left(\overline{\lambda}1_A|Z\right)\right)\left(E(\varphi_\lambda 1_B|Z) + E\left(\overline{\lambda}1_B|Z\right)\right) - E(\varphi_\lambda 1_A 1_B|Z) + E\left(\overline{\lambda}1_A 1_B|Z\right)}_{=:p'_\lambda} \overset{\text{a.s.}}{=} 0$$

Just like in the proof of Theorem 6.2.1, we can see that $\lambda \mapsto p'_\lambda(\omega)$ is a polynomial for all $\omega \in \Omega$. Since $\lambda \mapsto R_\lambda$ is continuous, there is $\varepsilon > 0$, s.t. $\forall \lambda < \varepsilon, R_\lambda \in \mathcal{D}$. Therefore, we have $p'_\lambda \overset{\text{a.s.}}{=} 0$ for all $\lambda < \varepsilon$. Therefore almost all polynomials $\lambda \mapsto p'_\lambda(\omega)$ are the zero polynomial. Therefore, $p'_1 \overset{\text{a.s.}}{=} 0$. Therefore, $A \perp\!\!\!\perp_{R_1} B \mid Z$. Since $R_1 = Q_{\leq n}$ and $A$ and $B$ were arbitrary, we have $X \perp\!\!\!\perp_{Q_{\leq n}} Y \mid Z$. □

**Corollary 6.4.4**: Suppose $X \not\perp Y \mid Z$. Then, in the topology induced by $d_1$, $\{P : X \perp\!\!\!\perp_P Y \mid Z\}$ is closed and nowhere dense and conversely, $\{P : X \not\perp\!\!\!\perp_P Y \mid Z\}$ is open and dense.

*Proof*: It suffices to show that $\{P : X \not\perp\!\!\!\perp_P Y \mid Z\}$ dense and $\{P : X \perp\!\!\!\perp_P Y \mid Z\}$ is closed. The former is a restatement of Theorem 6.4.3. The latter is Lemma 6.4.2. □

# 7 Properties of the history and structural independence

First, we take a look at the properties of structural independence. More precisely, the induced independence structure is a compositional semigraphoid.

In this section let $X, Y, Z, W$ be $\sigma(U)$-measurable random elements. Otherwise, the setting is taken from Section 6.

The goal of this section is to understand the properties of the history and prove that Desiderata 5.1 fully characterizes the history.

**Lemma 7.1** (monotonicity): If $\sigma(X,Z) \subseteq \sigma(Y,Z)$, then $\mathcal{H}(X|Z) \overset{\text{a.s.}}{\subseteq} \mathcal{H}(Y|Z)$.

*Proof*: It suffices to show that $\mathcal{H}(Y|Z)$ generates $X$ given $Z$. Firstly, by definition $\mathcal{H}(Y|Z)$ disintegrates $Z$, secondly, $\sigma(X,Z) \subseteq \sigma(Y,Z) \subseteq \sigma\left(\pi_{\mathcal{H}(Y|Z)}, Z\right)$. □

**Corollary 7.2** (monotonicity): if $\sigma(X) \subseteq \sigma(Y)$, then $\mathcal{H}(X|Z) \overset{\text{a.s.}}{\subseteq} \mathcal{H}(Y|Z)$.

*Proof*: This follows from Lemma 7.1. □

**Lemma 7.3** (compositionality): $\mathcal{H}((X,Y)|Z) \overset{\text{a.s.}}{=} \mathcal{H}(X|Z) \cup \mathcal{H}(Y|Z)$.

*Proof*: '$\overset{\text{a.s.}}{\supseteq}$' follows by Lemma 7.1.
'$\overset{\text{a.s.}}{\subseteq}$' Firstly, $H_1 = \mathcal{H}(X|Z)$ and $H_2 = \mathcal{H}(Y|Z)$ disintegrate $Z$ by definition. By the symmetry of disintegration, $H_1^c$ and $H_2^c$ disintegrate $Z$. By Lemma 5.5, we have $H_1^c \cap H_2^c$ disintegrates $Z$. Again, by symmetry, $H_1 \cup H_2$ disintegrate $Z$. Now $\sigma(X,Y) = \sigma(\sigma(X) \cup \sigma(Y)) \subseteq \sigma(\sigma\left(\left(U_{H_1}, Z\right) \cup \sigma\left(U_{H_2}, Z\right)\right) = \sigma\left(U_{H_1 \cup H_2}, Z\right)$ by Lemma 4.3.6.□

**Lemma 7.4** (idempotence): $\mathcal{H}\left(U_{\mathcal{H}(X|Z)}|Z\right) \overset{\text{a.s.}}{=} \mathcal{H}(X|Z)$.

*Proof*: Let $H := \mathcal{H}(X|Z)$.
'$\overset{\text{a.s.}}{\supseteq}$': By definition, $\sigma(X) \subseteq \sigma(U_H)$. By Lemma 7.1, $H \overset{\text{a.s.}}{\subseteq} \mathcal{H}(U_H|Z)$.
'$\overset{\text{a.s.}}{\subseteq}$': It suffices to show that $H$ generates $U_H$ given $Z$. By definition, $H$ disintegrates $Z$. Furthermore, $\sigma(U_H) \subseteq \sigma(U_H, Z)$. □

**Theorem 7.5**: Structural independence forms a compositional semigraphoid. More precisely, the following relations hold.

1. $X \perp Y \mid Z \Leftrightarrow Y \perp X \mid Z$ (symmetry)
2. $X \perp (Y, W) \mid Z \Rightarrow X \perp Y \mid Z$ (decomposition)
3. $X \perp (Y, W) \mid Z \Rightarrow X \perp Y \mid (Z, W)$ (weak union)
4. $X \perp Y \mid Z \wedge X \perp W \mid (Z, Y) \Rightarrow X \perp (Y, W) \mid Z$ (contraction)
5. $X \perp Y \mid Z \wedge X \perp W \mid Z \Rightarrow X \perp (Y, W) \mid Z$ (composition)

Here, 1-4. correspond to the semigraphoid axioms and 5 corresponds to the prefix 'compositional'.

*Proof*: It is well known that 1-4. hold for probability independence, see [9], Section 2 (6a)-(6e). By Theorem 6.4.1, we immediately obtain 1-4: We only prove 1. exemplary. It suffices to prove '$\Rightarrow$' Let $X \perp Y \mid Z$. Then $\forall P \in \triangle^{\times} : X \perp\!\!\!\perp_{\mathbb{P}} Y \mid Z$. Since 1. holds for $\perp\!\!\!\perp_{\mathbb{P}}$ in place of $\perp$, we have $\forall P \in \triangle^{\times} : Y \perp\!\!\!\perp_{\mathbb{P}} X \mid Z$ and therefore $Y \perp X \mid Z$. 2-4. are proved in the same way.

Finally, for 5. let $X \perp Y \mid Z$ and $X \perp W \mid Z$. Then $\mathcal{H}(X|Z) \cap \mathcal{H}(Y|Z) \overset{\text{a.s.}}{=} \emptyset$ and $\mathcal{H}(X|Z) \cap \mathcal{H}(W|Z) \overset{\text{a.s.}}{=} \emptyset$. Since $\mathcal{H}(Y, W|Z) \overset{\text{a.s.}}{=} \mathcal{H}(Y|Z) \cup \mathcal{H}(W|Z)$ by Lemma 7.3, it follow that $\mathcal{H}(X|Z) \cap \mathcal{H}(Y, W|Z) \overset{\text{a.s.}}{=} \emptyset$ and therefore $X \perp (Y, W) \mid Z$. □

Here it is important to point out that the composition axiom does not hold for probabilistic independence. This is because we will now see that structural independence does not differentiate between pairwise and 'full' independence of a vector.

**Definition 7.6** (structural independence of a vector): Let $(X_k)_{k \in K}$ be a family of $\sigma(U)$-measurable random elements. We say that $(X_k)_{k \in K}$ is structurally independent given $Z$, if for $k_1 \neq k_2 \in K$, we have $X_{k_1} \perp X_{k_2} \mid Z$.

**Theorem 7.7** (the fundamental theorem of structural independence for vectors): Let $X = (X_k)_{k \in K}$ be a family of $\sigma(U)$-measurable random elements. Then $X$ is structurally independent given $Z$ if and only if $X$ is independent given $Z$ for all distributions in $\triangle^{\times}$.

*Proof*: '$\Rightarrow$': W.l.o.g. $K$ is finite. Let $P \in \triangle^{\times}$. Then, by induction, we can assume $K = \{1, ..., n\}$ and $X' = (X_k)_{k=1}^{n-1}$ is independent. Let $A_k \in \sigma(X_k)$ for $k \in K$. Then by Theorem 7.5 and $X_n \perp X_k \mid Z$ for $k < n$, $X_n \perp X' \mid Z$. and therefore $P\left(\bigcap_{k \in K} A_k | Z\right) = P(A_n|Z) P\left(\bigcap_{k=1}^{n-1} A_k | Z\right) = \prod_{k \in K} P(A_k|Z)$.
'$\Leftarrow$': By decomposition of probabilistic independence, we have $\forall P \in \triangle^{\times} : X_{k_1} \perp\!\!\!\perp X_{k_2} \mid Z$ for all $k_1 \neq k_2 \in K$. By Theorem 6.4.1, $X_{k_1} \perp X_{k_2}$. □

We also have a lemma similar to how $X \perp\!\!\!\perp_{\mathbb{P}} Y \mid Z$ implies $\mathbb{E}(X|Z) \overset{\text{a.s.}}{=} \mathbb{E}(X|Y, Z)$.

**Lemma 7.8**: If $X \perp Y \mid Z$, then $\mathcal{H}(X|Z) \overset{\text{a.s.}}{=} \mathcal{H}(X|Y, Z)$.

*Proof*: Clearly, $H := \mathcal{H}(X|Z)$ is a $\sigma(Y, Z)$-measurable random index set. Furthermore, $\sigma(X) \subseteq \sigma(U_H, Z) \subseteq \sigma(U_H, Y, Z)$. Therefore, it suffices to show that $H$ disintegrates $(Y, Z)$. Now, by definition, we have $U_H \perp U_{\overline{H}} \mid Z$. By Lemma 7.4, we also have $U_H \perp Y \mid Z$. Therefore, $U_H \perp (Y, U_{\overline{H}}) \mid Z$. By Theorem 7.5 (weak union), we have $U_H \perp U_{\overline{H}} \mid (Y, Z)$. □

**Lemma 7.9**: We have $\mathcal{H}(X|Y, Z) \overset{\text{a.s.}}{\subseteq} \mathcal{H}(X, Y|Z)$.

*Proof*: It suffices to show that $H := \mathcal{H}(X, Y|Z)$ generates $X$ given $(Y, Z)$. By definition, we have $\sigma(X) \subseteq \sigma(X, Y) \subseteq \sigma(U_H, Z) \subseteq \sigma(U_H, Y, Z)$. It is left to show that $H$ disintegrates $(Y, Z)$. By definition, we have $H$ disintegrates $Z$. Therefore $U_H \perp U_{\overline{H}} \mid Z$. Since $\mathcal{H}(Y|Z) \cap \overline{H} \overset{\text{a.s.}}{=} \emptyset$ and, since $\overline{H}$ disintegrates $Z$, $\mathcal{H}(U_{\overline{H}}) \subseteq \overline{H}$, we have $Y \perp U_{\overline{H}} \mid Z$. Therefore, $U_{\overline{H}} \perp (Y, U_H) \mid Z$. By Theorem 7.5 (weak union), we have $U_{\overline{H}} \perp U_H \mid (Y, Z)$. □

We will now prove that Desiderata 5.1 uniquely determines the history. For this we first show that the history fulfills these desiderata.

**In the following**, let $J$ be a random index set.

**Lemma 7.10**: Let $J$ disintegrate $Z$ and $K = \bigcup^{\text{a.s.}}\{K' : K'$ is a $\sigma(Z)$-measurable random index set s.t. $\sigma(U_{K'}) \subseteq \sigma(Z)\}$. Then $\mathcal{H}(U_J|Z) = J \setminus K$.

*Proof*: '$\overset{\text{a.s.}}{\subseteq}$': By Lemma 4.3.6, $\sigma(U_K) \subseteq \sigma(Z)$. Then clearly $K$ disintegrates $Z$. By Lemma 5.5 and the symmetry of disintegration, $J \setminus K$ disintegrates $Z$. Clearly, $\sigma(U_J, Z) = \sigma(U_{J \setminus K}, Z)$.
'$\overset{\text{a.s.}}{\supseteq}$': Let $L = (J \setminus K) \setminus \mathcal{H}(U_J|Z)$. By Lemma 5.5 and the symmetry of disintegration, $L$ disintegrates $Z$. By '$\overset{\text{a.s.}}{\subseteq}$', we have $\mathcal{H}(U_L|Z) \overset{\text{a.s.}}{\subseteq} L$. By Theorem 6.4.1, we have $U_L \perp\!\!\!\perp U_J \mid Z$ and therefore $U_L \perp\!\!\!\perp U_L \mid Z$, therefore $\sigma(U_L) \subseteq \sigma(Z)$. Therefore, $L \overset{\text{a.s.}}{\subseteq} K$. But by definition, $L \cap K \overset{\text{a.s.}}{=} \emptyset$, therefore $L \overset{\text{a.s.}}{=} \emptyset$. □

**Corollary 7.11**:

$$\overline{\mathcal{H}(U|Z)} \overset{\text{a.s.}}{=} \bigcup^{\text{a.s.}}\{K' : K' \text{ is a } \sigma(Z)\text{-measurable random index set s.t. } \sigma(U_{K'}) \subseteq \sigma(Z)\}.$$

*Proof*: Follows immediately from Lemma 7.10. □

**Lemma 7.12**: $\mathcal{H}(X|Z) \overset{\text{a.s.}}{=} \emptyset \Leftrightarrow \sigma(X) \subseteq \sigma(Z)$.

*Proof*: Follows immediately from the definition of the history. □

**Lemma 7.13**: $\mathcal{H}\left(U_{J \setminus \mathcal{H}(U_J|Z)}|Z\right) \overset{\text{a.s.}}{=} \emptyset$

*Proof*: Let $K = J \setminus \mathcal{H}(U_J|Z)$. By Lemma 7.1 and Lemma 7.10, we have $\mathcal{H}(U_K|Z) \overset{\text{a.s.}}{\subseteq} \mathcal{H}\left(U_{\overline{\mathcal{H}(U_J|Z)}}|Z\right) \overset{\text{a.s.}}{\subseteq} \overline{\mathcal{H}(U_J|Z)}$. By Theorem 6.4.1, $U_K \perp U_J \mid Z$. Since $K \subseteq J$, we have $U_K \perp U_K \mid Z$, therefore $\sigma(U_K) \subseteq \sigma(Z)$. By Lemma 7.12, $\mathcal{H}(U_K) \overset{\text{a.s.}}{=} \emptyset$. □

**Lemma 7.14**: If $\mathcal{H}(U_J|Z) \overset{\text{a.s.}}{=} \emptyset$, then $\mathcal{H}(U|Z) \cap J \overset{\text{a.s.}}{=} \emptyset$.

*Proof*: By Lemma 7.12, $\sigma(U_J) \subseteq \sigma(Z)$. Then it is clear that $\overline{J}$ generates $U$ given $Z$. □

**Lemma 7.15**: $\mathcal{H}(U|Z) \cap J \overset{\text{a.s.}}{\subseteq} \mathcal{H}(U_J|Z)$.

*Proof*: By Lemma 7.13 and Lemma 7.14, we have $\mathcal{H}(X|Z) \cap J \setminus \mathcal{H}(U_J) \overset{\text{a.s.}}{=} \emptyset$. Setting $X = U$, we have $\mathcal{H}(U|Z) \cap J \overset{\text{a.s.}}{\subseteq} \mathcal{H}(U_J|Z)$. □

We can now prove Desiderata 5.1.

**Theorem 7.16** (maximality and uniqueness of the history): Let $\Sigma$ be the set of complete sub-$\sigma$-algebras of $\mathcal{A}$. The history is the almost surely maximal and unique map $\mathcal{H}(\cdot \mid \cdot) : \Sigma \times \Sigma \to \mathfrak{P}(I)^\Omega$ that fulfills the following,

1. $\mathcal{H}(X|Z)$ is a $\sigma(Z)$-measurable random index set.
2. $\forall P \in \triangle^\times : X \perp\!\!\!\perp_P Y \mid Z \Leftrightarrow \mathcal{H}(X|Z) \cap \mathcal{H}(Y|Z) \overset{\text{a.s.}}{=} \emptyset$.
3. If $\sigma(X) \subseteq \sigma(Y)$ then $\mathcal{H}(X|Z) \overset{\text{a.s.}}{\subseteq} \mathcal{H}(Y|Z)$.
4. $\mathcal{H}(U|Z) \cap J \overset{\text{a.s.}}{\subseteq} \mathcal{H}(U_J|Z)$ for all random index sets $J$.

*Proof*: The history fulfills 1. by definition, 2. by Theorem 6.4.1, 3., by Lemma 7.1 and 4. by Lemma 7.15.

'Maximality'. Let $\mathcal{H}'$ be a map that fulfills 1-4. Let $Z$ be fixed. We need to show that $\mathcal{H}'(X|Z) \overset{\text{a.s.}}{\subseteq} \mathcal{H}(X|Z)$ for all $X$.

i) Let $J$ disintegrate $Z$. We show that $\mathcal{H}'(U_J, Z|Z) \overset{\text{a.s.}}{\subseteq} J$. By 4. applied to $\overline{J}$, we have
$$\mathcal{H}'(U|Z) \cap \overline{J} \overset{\text{a.s.}}{\subseteq} \mathcal{H}'(U_{\overline{J}}|Z)$$
$$\Leftrightarrow \mathcal{H}'(U|Z) \cap \overline{J} \cap \overline{\mathcal{H}'(U_{\overline{J}}|Z)} \overset{\text{a.s.}}{=} \emptyset$$
Since $(U_J, Z) \perp U_{\overline{J}} \mid Z$, we have $\mathcal{H}'(U_J, Z|Z) \overset{\text{a.s.}}{\subseteq} \overline{\mathcal{H}'(U_{\overline{J}}|Z)}$ by 2. Therefore we have by 3. and the $\sigma(U)$-measurability of $Z$, that
$$\begin{aligned} \mathcal{H}'(U, Z|Z) \cap \overline{J} \cap \mathcal{H}'(U_J, Z|Z) &\overset{\text{a.s.}}{=} \emptyset \\ \Leftrightarrow \qquad \mathcal{H}'(U_J, Z|Z) \cap \overline{J} &\overset{\text{a.s.}}{=} \emptyset \\ \Leftrightarrow \qquad \mathcal{H}'(U_J, Z|Z) &\overset{\text{a.s.}}{\subseteq} J \end{aligned}$$

ii) We show that $\mathcal{H}'(Z|Z) \overset{\text{a.s.}}{=} \emptyset$.
Clearly, $Z \perp Z \mid Z$. By 2. $\mathcal{H}'(Z|Z) = \mathcal{H}'(Z|Z) \cap \mathcal{H}'(Z|Z) \overset{\text{a.s.}}{=} \emptyset$.

iii) We show that $\mathcal{H}'(U|Z) \overset{\text{a.s.}}{\subseteq} \mathcal{H}(U|Z)$. By Corollary 7.11,
$$\overline{\mathcal{H}(U|Z)} = K := \overset{\text{a.s.}}{\bigcup} \{K' : K' \text{ is a } \sigma(Z)\text{-measurable random index set s.t. } \sigma(U_{K'}) \subseteq \sigma(Z)\}.$$
It suffices to show $L := K \cap \mathcal{H}'(U|Z) \overset{\text{a.s.}}{=} \emptyset$. By 4. and 3. we have $L \overset{\text{a.s.}}{\subseteq} \mathcal{H}'(U_L|Z) \overset{\text{a.s.}}{\subseteq} \mathcal{H}'(U_K|Z)$. It suffices to show that $\mathcal{H}'(U_K|Z) \overset{\text{a.s.}}{=} \emptyset$. By definition of $K$, $\sigma(U_K) \overset{\text{a.s.}}{\subseteq} \sigma(Z)$. By 3. and ii) we have $\mathcal{H}'(U_K|Z) \overset{\text{a.s.}}{\subseteq} \mathcal{H}'(Z|Z) \overset{\text{a.s.}}{=} \emptyset$.

iv) We show that $\mathcal{H}'(U_J, Z|Z) \overset{\text{a.s.}}{=} J \cap \mathcal{H}'(U|Z)$, whenever $J$ disintegrates $Z$.
'$\overset{\text{a.s.}}{\subseteq}$': By 3. $\mathcal{H}'(U_J, Z|Z) \overset{\text{a.s.}}{\subseteq} \mathcal{H}'(U|Z)$. By i) we have $\mathcal{H}'(U_J, Z|Z) \overset{\text{a.s.}}{\subseteq} J$.
'$\overset{\text{a.s.}}{\supseteq}$': By 4. we have $\mathcal{H}'(U|Z) \cap J \overset{\text{a.s.}}{\subseteq} \mathcal{H}'(U_J|Z) \overset{\text{a.s.}}{\subseteq} \mathcal{H}'(U_J, Z|Z)$.

v) We claim that $\mathcal{H}'(U_J, Z|Z) \overset{\text{a.s.}}{\subseteq} \mathcal{H}(U_J, Z|Z)$ whenever $J$ disintegrates $Z$. Note that iii) shows $\mathcal{H}'(U|Z) \overset{\text{a.s.}}{\subseteq} \mathcal{H}(U|Z)$. Now $\mathcal{H}'(U_J, Z|Z) = J \cap \mathcal{H}'(U|Z) \overset{\text{a.s.}}{\subseteq} J \cap \mathcal{H}(U|Z) = \mathcal{H}(U_J, Z|Z)$.

vi) Finally, let $X$ be arbitrary. We show that $\mathcal{H}'(X|Z) \overset{\text{a.s.}}{\subseteq} \mathcal{H}(X|Z)$. Let $H := \mathcal{H}(X|Z)$
Since $X \overset{\text{a.s.}}{\subseteq} \sigma(U_H, Z)$ by definition of $H$, by 3. we have $\mathcal{H}'(X|Z) \overset{\text{a.s.}}{\subseteq} \mathcal{H}'(U_H, Z|Z)$. By v), $\mathcal{H}'(U_H, Z|Z) \overset{\text{a.s.}}{=} H \cap \mathcal{H}'(U|Z) \overset{\text{a.s.}}{\subseteq} H$.

'Uniqueness': Let $\mathcal{H}^1, \mathcal{H}^2$ be two maximal maps that fulfill 1-4. Then by maximality, we have $\mathcal{H}^1(X|Z) \overset{\text{a.s.}}{\subseteq} \mathcal{H}^2(X|Z) \overset{\text{a.s.}}{\subseteq} \mathcal{H}^1(X|Z)$. □

# 8 A counterexample

In this section, we introduce a example that shows that disintegrates cannot be characterized by rectangular atoms in general, in contrast to the finite case, see Section 3. Furthermore, it illustrates that the choice of a reference measure (or an equivalence class of mutually absolutely continuous probability measures) is necessary for the history to exist. For this we introduce a certain $Z$ on a two-dimensional product space.

**Example 8.1**: Let $S = [0, 1]$, Let $I = \{1, 2\}$ and $i \in I$. Let $\Omega_i = S \sqcup S$, where $\sqcup$ denotes the disjoint union, i.e. $A \sqcup B = A \times \{1\} \cup B \times \{2\}$. To access the two parts of $\Omega_i$, we write $S_i = S \times \{i\}$. Let $\mathcal{A}_i$ be the Borel $\sigma$-algebra on $\Omega_1$. Set $(\Omega, \mathcal{A}) = (\Omega_1 \times \Omega_2, \mathcal{A}_1 \otimes \mathcal{A}_2)$. We set $U_i = \pi_i : \Omega \to \Omega_i$ to be the canonical projection.

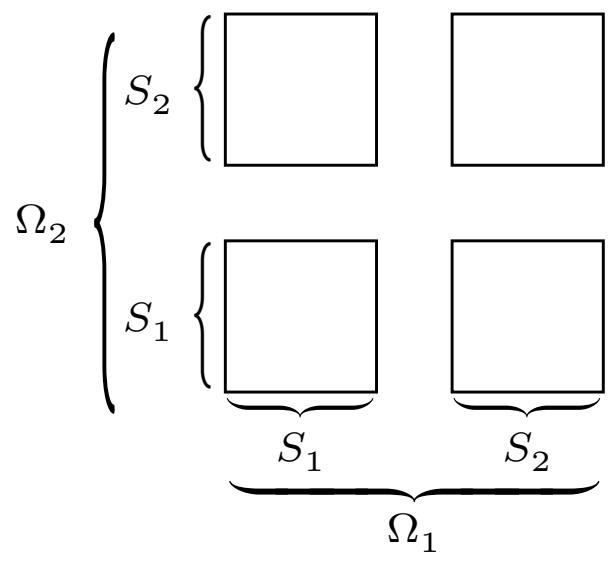

Figure 1: An illustration of $\Omega = \Omega_1 \times \Omega_2$

We now construct the random element $Z : \Omega \to S^2$ on which we will condition. For $i, j \in I$, set $S_{ij} = S_i \times S_j$. Then $\Omega = \bigcup_{i,j \in I} S_{ij}$. Therefore it suffices to define $Z$ on each of $S_{ij}$, we write $Z_{ij}$ for $Z|_{S_{ij}}$.

Let $\alpha \in (0,1)$ and $\beta \in (0,1)$. Let $E = \{(a,b) \in S^2 : a + \beta \cdot b, \alpha\alpha + \beta \in [0,1]\}$. Let for $(a,b) \in E$

$$Z_{12}^{-1}\begin{pmatrix} a \\ b \end{pmatrix} = \begin{pmatrix} 1 & 0 \\ \alpha & 1 \end{pmatrix}\begin{pmatrix} a \\ b \end{pmatrix} = \begin{pmatrix} a \\ \alpha \cdot a + b \end{pmatrix} \in S_{12} \qquad Z_{22}^{-1}\begin{pmatrix} a \\ b \end{pmatrix} = \begin{pmatrix} 1 & \alpha \\ \beta & 1 \end{pmatrix}\begin{pmatrix} a \\ b \end{pmatrix} = \begin{pmatrix} a + \beta \cdot b \\ \alpha \cdot a + b \end{pmatrix} \in S_{22}$$

$$Z_{11}^{-1}\begin{pmatrix} a \\ b \end{pmatrix} = \begin{pmatrix} 1 & 0 \\ 0 & 1 \end{pmatrix}\begin{pmatrix} a \\ b \end{pmatrix} = \begin{pmatrix} a \\ b \end{pmatrix} \in S_{11} \qquad Z_{21}^{-1}\begin{pmatrix} a \\ b \end{pmatrix} = \begin{pmatrix} 1 & \beta \\ 0 & 1 \end{pmatrix}\begin{pmatrix} a \\ b \end{pmatrix} = \begin{pmatrix} a + \beta \cdot b \\ b \end{pmatrix} \in S_{21}$$

The partial functions $Z_{ij} : S_{ij} \rightharpoonup E$ are well defined as the inverses of $Z_{i,j}^{-1}$, since $Z_{ij}^{-1}$ is bijective by $\det\begin{pmatrix} 1 & \beta \\ \alpha & 1 \end{pmatrix} = 1 - \alpha\beta \neq 0$, etc. Let $D = \bigcup_{i,j \in I} Z_{ij}^{-1}(E))$ denote the union of their support. Since, $Z_{ij}^{-1}$ have disjoint support, we can define $Z : \Omega \to S^2 \sqcup \Omega$ by $Z|_D = \bigcup_{ij} Z_{ij} : D \to S^2$ and $Z|_{D^c} = \mathrm{id}|_{D^c} : D^c \to \Omega$. Therefore $\sigma(Z|_{D^c}) = \sigma(\mathcal{A}|_{D^c})$, so that we only have to consider $Z|_D$ when checking for conditional independence of $U_1$ and $U_2$ given $Z$.

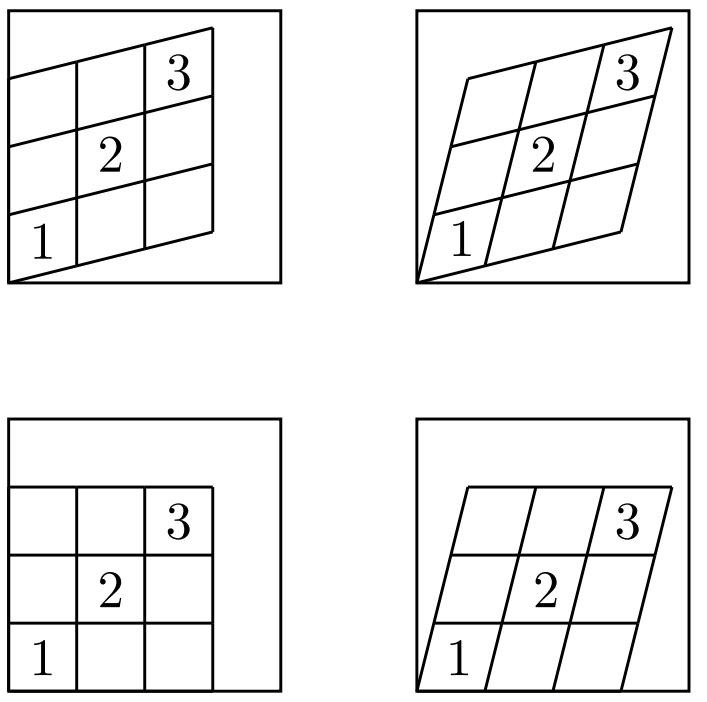

Figure 2: An illustration of $Z$ for $\alpha = \beta = \frac{1}{4}$ in the region $Z^{-1}\left[0, \frac{3}{4}\right]^2$. The cells correspond to the partition $\left\{Z \in [a - \varepsilon, a + \varepsilon) \times [b - \varepsilon, b + \varepsilon) : a, b \in \left\{\varepsilon\left(n + \frac{1}{2}\right) : n \in \{0, 1, 2\}\right\}\right\}$ for $\varepsilon = \frac{1}{4}$. The numbers inscribed in the cells illustrate which cells are in the same part of the partition.

Note that $Z$ has rectangular atoms, since for $(a,b) \in E$,

$$Z^{-1}(a,b) = \{a, a + \beta \cdot b\} \times \{b, \alpha \cdot a + b\}$$

and otherwise, $Z^{-1}(a,b)$ is a single element set or empty.

Now let $\mathbb{P} = \mathbb{P}_1 \times \mathbb{P}_2$, where $\mathbb{P}_i$ is the uniform distribution over $S \sqcup S$. I.e. $\mathbb{P}|_{S_{ij}} = \frac{1}{4}\lambda$, where $\lambda$ is the Lebesgue measure over $S_{ij}$. We claim that the constant random index set $\{1\}$ does not disintegrate $Z$.

Let $A = S_1 \times \Omega_2$ and $B = \Omega_1 \times S_1$. Then $A \cap B = S_{11}$. It suffices to show that $\mathbb{P}(A|Z)\mathbb{P}(B|Z) \overset{\text{a.s.}}{\neq} \mathbb{P}(S_{11}|Z)$. We use the Lebesgue differentiation theorem to calculate the conditional expectations, cf [10], Section 7.2. For $(a, b) \in S^2, \varepsilon \in \mathbb{R}$, let $B_\varepsilon(a, b) = (a - \varepsilon, a + \varepsilon) \times (b - \varepsilon, b + \varepsilon)$ denote the $\varepsilon$ ball around $(a, b)$ in $S^2$. For $z \in S^2$, let $B_\varepsilon^z = \{Z \in B_\varepsilon(z))\}$.

In Figure 2, we can already see that the volume the partition corresponding to 1, as $\varepsilon \to 0$ is not of a product structure: Let $P$ denote the part in the partition corresponding to 1. Then, in the bottom left ($P \cap S_{11}$), the volume is 1, and so is the top left ($P \cap S_{12}$) and bottom right ($P \cap S_{21}$). On the other hand, the top right ($P \cap S_{22}$) is of a diamond shape, so that the volume is not 1. Therefore $A$ and $B$ will not be independent given $Z$. We make this more precise now.

Let $L$ be a linear invertible map and $B$ a Borel set. Let $\lambda$ be the Lebesgue measure. Then $\lambda(L(B)) = \lambda(B) \cdot |\det(L)|$ Using this formula, and noting that $\det Z_{ij}^{-1} = 1$ for $(i, j) \neq (2, 2)$ and $\det Z_{22}^{-1} = (1 - \alpha\beta)$, we have, for $z \in (0, 1)^2$ and $\varepsilon$ small enough, that $\mathbb{P}\big(S_{ij} \cap B_\varepsilon^z\big) = \mathbb{P}\big(Z_{ij}^{-1}(B_\varepsilon(z))\big) = \frac{1}{4}\lambda(B_\varepsilon(z)) = \frac{\varepsilon^2}{4}$ for $(i, j) \neq (2, 2)$ and $\mathbb{P}(S_{22} \cap B_\varepsilon^z) = \mathbb{P}\big(Z_{22}^{-1}(B_\varepsilon(z))\big) = \frac{\varepsilon^2}{4}(1 - \alpha\beta)$.

Therefore the following holds for almost all $z \in D$

$$\mathbb{P}(A|Z = z) = \lim_{\varepsilon \to 0} \mathbb{P}(A \cap B_\varepsilon^z)/\mathbb{P}(B_\varepsilon^z) = \lim_{\varepsilon \to 0} \frac{\frac{2}{4}\varepsilon^2}{\big(\frac{3}{4} + \frac{1}{4}(1 - \alpha\beta)\big)\varepsilon^2} = \frac{2}{4 - \alpha\beta}.$$

Similarly, $\mathbb{P}(B|Z = z) = \frac{2}{4-\alpha\beta}$ and $\mathbb{P}(A, B|Z = z) = \frac{1}{4-\alpha\beta}$. Now, setting $x = \alpha\beta$

$$\begin{aligned} \mathbb{P}(A|Z = z)\mathbb{P}(B|Z = z) = \mathbb{P}(A, B|Z = z) &\Leftrightarrow \left(\frac{2}{4 - x}\right)^2 = \frac{1}{4 - x} \\ &\Leftrightarrow 4 = 4 - x \end{aligned}$$

This is clearly false, since $x = \alpha\beta > 0$. Therefore, $\mathbb{P}(A|Z)\mathbb{P}(B|Z) \neq P(A, B|Z)$ a.s. on $D$.

In conclusion, we have seen that there is a $Z$, s.t. $Z^{-1}(z)$ is a rectangle for all $z$ in the codomain of $Z$, but also $\{1\}$ does not disintegrate $Z$.

**Example 8.2**: We continue the previous example. Let $\Sigma$ be the set of sub-$\sigma$-algebras of $\mathcal{A}$. Notably, here we don't complete them w.r.t. a reference measure. Let $\mathcal{H}_P(X|Z)$ denote the history, defined through the reference measure $P$. We use the example to show that we there is no map $\mathcal{H}(\cdot, \cdot) : \Sigma \times \Sigma \to \mathfrak{P}(I)^\Omega$, s.t. $\mathcal{H}(X|Z) = \mathcal{H}_P(X|Z)$ $P$-a.s. for all $P$ s.t. $U$ is independent w.r.t. $P$. Suppose there is such a map $\mathcal{H}$.

Recall that $\mathbb{P}$ is the uniform distribution. Clearly, $\sigma(U_1|_C) \not\subseteq \sigma(Z)$ for any $\mathbb{P}$-non-nullset $C \subseteq D$. Using the same arguments as in the previous example, we can see that the only disintegrating random index sets (w.r.t. $\mathbb{P}$) are $\emptyset$ and $I$. $\mathbb{P}$-a.e. on $D$. Therefore, $\mathcal{H}(U_1|Z)(\omega) = I$ for $\mathbb{P}$-a.e. $\omega \in D$

We will now show that $\mathcal{H}(U_1|Z) = \emptyset$ everywhere, contradicting the existence of $\mathcal{H}$ by the previous paragraph. Let $z \in \mathrm{Val}(Z)$ and let $C = Z^{-1}\{z\}$. Since $C$ is a rectangle with four elements, $P = \frac{1}{4}\sum_{c \in C} \delta_c$, where $\delta_c$ is the Dirac measure, is a product probability measure on $(\Omega, \mathcal{A})$. By Definition 3.1, $\{1\}$ disintegrates $Z$ (w.r.t. $P$). Therefore, clearly, $\mathcal{H}_P(U_1|Z) = 1$ $P$-a.s. Now $\mathcal{H}(U_1|Z) = \mathcal{H}_P(U_1|Z)$ $P$-a.s. and therefore $\mathcal{H}(U_1|Z)(\omega) = \{1\}$ for $\omega \in C$. But since $z$ was arbitrary and $\{Z^{-1}\{z\} : z \in \mathrm{Val}(z)\}$ covers $\Omega$, we have $\mathcal{H}(U_1|Z)(\omega) = \{1\}$ for all $\omega \in \Omega$. But then clearly, $\mathcal{H}(U_1|Z) \neq \mathcal{H}_{\mathbb{P}}(U_1|Z)$ $\mathbb{P}$-a.s. This is a contradiction to the existence of $\mathcal{H}$.

# 9 Applications

In this section we discuss the implications of structural independence for structural causal models. The following definition is equivalent to [2], Definition 7.1.1, Definition 7.1.6

**Definition 9.1** ((Structural) causal model): A (structural) causal model (SCM) is a triple $(U, V, F)$, where

1. $U = (U_i)_{i \in I}$ is a family of random variables. These variables are called exogenous or background variables.
2. $V = (V_i)_{i \in I}$ is a family of random variables, s.t. $\sigma(V) \subseteq \sigma(U)$. These variables are called endogenous or observed variables.
3. $F = (F_i)_{i \in I}$ is a family of functions $F_i : \mathrm{Val}(U_i) \times \mathrm{Val}\left(V_{\mathrm{PA}(i)}\right) \to \mathrm{Val}(V_i)$, where $\mathrm{PA}(i) \subseteq I$.

All random variables have common domain $(\Omega, \mathcal{A})$. and satisfy

$$\forall i \in I : V_i = F_i\left(U_i, V_{\mathrm{PA}(i)}\right)$$

If $\mathbb{P}$ is a probability distribution on $(\Omega, \mathcal{A})$, s.t. $U$ is an independent family w.r.t. $\mathbb{P}$, we call $(U, V, F, \mathbb{P})$ a probabilistic causal model.

To a causal model $\mathcal{M}$ we associate a graph $G_{\mathcal{M}} = (V, E)$, where $E = \left\{\left(V_j, V_i\right) \in V \times V : j \in \mathrm{PA}(i)\right\}$.

An SCM $\mathcal{M}$ is acyclic (or recursive) if $G_{\mathcal{M}}$ is acyclic. In this paper we always assume that the causal models are acyclic.

We call an acyclic probabilistic causal model purely probabilistic, if $V_i \not\!\perp\!\!\!\perp U_i$ for all $i$. I.e. $U_i$ is needed to determine $V_i$. Otherwise, $V_i$ would be a deterministic function of $V_{\mathrm{PA}(V_i)}$.

We interpret a causal model as follows: $U$ are the background variables that are not observed. $V$ are the observed variables. $F_i$ represents the mechanism that determines $V_i$ noisily. The noise is given by $U_i$. $F_i$ embodies that $V_{\mathrm{PA}(V_i)}$ cause $V_i$, since $V_i$ is a probabilistic function of $V_{\mathrm{PA}(V_i)}$.

We note that a (causal) bayesnet can be understood as a special case of a causal model. More specifically, we can define interventions on causal bayesnets and on structural causal models. We can construct a causal model from a graph s.t. these interventions agree. For more details, see [8], [11, Section 4].

It is well known how to do simple causal discovery in the case of graphs from the combinatorial properties of a given probability distribution (c.f. Section 2 or [2, Chapter 2]), but the same type theory has not been developed for structural causal models. This is because until now, there was no equivalent to $d$-separation for structural causal models. Structural independence is this equivalent notion. We now state the fundamental theorem in the structural causal model setting.

**Theorem 9.2** (Fundamental theorem of structural independence in causal models): Let $(U, V, F, \mathbb{P})$ be a probabilistic structural causal model. Let $\mathcal{M} = \{(U, V, F, P) \text{ probabilistic causal model} : P \sim \mathbb{P}\}$ be the set of all probabilistic causal models that agree with $(U, V, F)$ and where $P$ is comparable to $\mathbb{P}$. Then for all random variables $X, Y, Z$ that are determined by $U, V$, i.e. $\sigma(X) \subseteq \sigma(U)$, we have

$$\forall (U, V, F, P) \in \mathcal{M} : X \perp\!\!\!\perp_P Y \mid Z \Leftrightarrow X \perp_U Y \mid Z$$

where $X \perp Y \mid Z :\Leftrightarrow \mathcal{H}(X|Z) \cap \mathcal{H}(Y|Z) \overset{\mathbb{P}\text{–a.s.}}{=} \emptyset$. $X \perp Y \mid Z$ is structural independence w.r.t. the independent family $U$ (Definition 5.12).

*Proof*: This is a restatement of Theorem 6.4.1. □

We are confident that this theorem can be leveraged for causal discovery in structural causal models, just like $d$-separation is used for causal discovery in bayesnets.

The following example is adapted from [6, Example 1] and [8] and shows that, in principle, structural independence can be used to discover causal relations in structural causal models from probabilistic data.

**Example 9.3**: We observe two nonconstant, (real valued) random variables $V_1$ and $V_2$. We assume that they are part of an unobserved true underlying causal model $(U, V, F, \mathbb{P})$, where $V = (V_i)_{i \in I}$ and $\{1, 2\} \subseteq I$. We only observe the distribution of $V$. Let $Z \coloneqq V_1 + V_2$. Suppose that $V_1 \perp\!\!\!\perp_{\mathbb{P}} Z$. Furthermore, we assume that $\mathbb{P}$ is in 'general position'. This means that this independence is not due to the specific choice of $\mathbb{P}$. Formally, $V_1 \perp\!\!\!\perp_P Z$ for all causal models $(U, V, F, P)$ s.t. $P \sim \mathbb{P}$.

Then $V_1$ is an ancestor of $V_2$ in $G_{\mathcal{M}}$.

Furthermore, if $V = (V_1, V_2)$, then $V_1$ is a parent of $V_2$ in $G_{\mathcal{M}}$.

We use two lemmas to proof this.

**Lemma 9.4**: Let $\mathcal{M} = (U, V, F, \mathbb{P})$ be a purely probabilistic causal model. Let $\mathcal{H}(V_i) \overset{\text{a.s.}}{\subseteq} \mathcal{H}(V_j)$, then $V_i$ is an ancestor of $V_j$ in $G_{\mathcal{M}}$.

*Proof*: Since $\mathcal{M}$ is purely probabilistic, we have $i = \mathcal{H}(U_i) \subseteq \mathcal{H}(V_i) \subseteq \mathcal{H}(V_j)$ Now it is easy to see that since $U_i$ is needed to determine $V_j$ and $V_j$ can only access $U_i$ through $V_i$, $V_i$ needs to be an ancestor of $V_j$ in $G_{\mathcal{M}}$. □

We now state a more general version of Example 9.3

**Lemma 9.5**: Let $(G, *, \mathcal{G})$ be a measurable group. I.e. $\mathcal{G}$ is a $\sigma$-algebra over $G$, and the group operation $*$ and inversion $\cdot^{-1}$ is measurable. Let $(\Omega, \mathcal{A}, \mathbb{P})$ be a probability space. Let $V_1, V_2$ be $G$-valued random elements on $(\Omega, \mathcal{A}, \mathbb{P})$. Suppose that $V_1 \perp (V_1 * V_2)$, then $\mathcal{H}(V_1) \subseteq \mathcal{H}(V_2)$.

*Proof*: Since $V_1 \perp Z$, we have $\mathcal{H}(V_1) \cap \mathcal{H}(V_1 * V_2) = \emptyset$. Since $V_1 = (V_1 * V_2) * V_2^{-1}$, we have $\mathcal{H}(V_1) \subseteq \mathcal{H}(V_1 * V_2, V_2) = \mathcal{H}(V_1 * V_2) \cup \mathcal{H}(V_2)$. Since $\mathcal{H}(V_1) \cap \mathcal{H}(V_1 * V_2) = \emptyset$, we must have $\mathcal{H}(V_1) \subseteq \mathcal{H}(V_2)$. □

*Proof of Example 9.3*: To prove that $V_1$ is an ancestor of $V_2$, using Lemma 9.4, it suffices to show that $\mathcal{H}(V_1) \subseteq \mathcal{H}(V_2)$. This is precisely Lemma 9.5. If $V = (V_1, V_2)$, then if $V_1$ is an ancestor of $V_2$ in $G_{\mathcal{M}}$, we must have $V_1 \in \mathrm{PA}(V_2)$. □

# 10 Further work

In this section we discuss further work that can be done on the theory of structural independence. We continue the setting of Section 6.

## 10.1 Disintegration

In Lemma 5.3, we have seen that it suffices to check the disintegration condition for one $P \in \triangle^{\times}$. Therefore this criterion is testable once we have any probability distribution for which we want to test which independences are structural. Nonetheless, in the finite case, Section 3, there is a elegant characterization of disintegration, namely that atoms of the $\sigma$-algebra of the conditional $Z$ are rectangles w.r.t. the random index set $J$, i.e. $U(C) = U_J(C) \times U_{\overline{J}}(C)$ for all atoms $C$ of $\sigma(Z)$. It is possible that a similar characterization is possible for the infinite setting. We provide a necessary but not sufficient condition.

**Lemma 10.1.1**: If $J$ disintegrates $Z$, then for all $A \in \sigma(U_J, Z)$ and $B \in \sigma(U_{\overline{J}}, Z)$, s.t. $A \cap B \overset{\text{a.s.}}{=} \emptyset$ , there is a $C \in \sigma(Z)$, s.t. $A \subseteq C$ and $B \subseteq C^c$.

*Proof*: We have $\mathbb{P}(A|Z)\mathbb{P}(B|Z) = \mathbb{P}(A, B|Z) \overset{\text{a.s.}}{=} 0$. Let $C = \{\mathbb{P}(A|Z) > 0\}$. Then $A \overset{\text{a.s.}}{\subseteq} C$. Furthermore, $\mathbb{P}(B|Z)(\omega) = 0$ for a.e. $\omega \in C$. Therefore, $B \overset{\text{a.s.}}{\subseteq} \{\mathbb{P}(B|Z) > 0\} \overset{\text{a.s.}}{\subseteq} C^c$. □

It can be seen that in the finite case, this fully characterizes disintegration, since it encodes the rectangle condition mentioned before. However, in the general case, it does not. It can be seen that in Example 8.1, the condition holds for $J = \{1\}$, but $\{1\}$ does not disintegrate $Z$. The condition can be morally understood as every interaction between $U_J$ and $U_{\overline{J}}$ being mediated by $Z$. We contend that it fails to capture disintegration since in the infinite theory, a conditional probability can be understood as a limit procedure, while the condition only talks about 'stationary' sets.

It could be that we need to introduce limiting objects to capture these phenomena.

**Remark 10.1.2**: Let $(A_n)_{n\in\mathbb{N}}$ be a sequence in $\mathcal{A}$. Let $A \in \mathcal{A}$. We define $A_n \to A$ if $1_{A_n} \to 1_A$ $\mathbb{P}$-a.s. Then this clearly does not dependent anything but the nullsets of $\mathbb{P}$. In other words, this is a sense of convergence w.r.t. nullsets.

**Conjecture 10.1.3**: Let $\Omega$ be a polish space. There is a suitable sense of convergence of $\sigma$ -algebras that only depends on the nullsets of $\mathbb{P}$, s.t. $J$ disintegrates $Z$, if and only if there is a sequence of finitely generated $\sigma$-algebras $\mathcal{A}_n$ that have rectangular atoms and converges to $\sigma(Z)$.

## 10.2 Conditional systems

In this work, we have focused on the independences that are implied by a family of random element $U$ being independent, since this reflects how structural causal models are defined. It seems highly likely that the theory can be extended to a family that fulfills certain conditional independences instead. More precisely, let $I = \mathbb{N}$. For $i \in I$, choose a set $J_i \subseteq \{j \in I : j < i\}$. Let $\mathbb{P}$ be a reference measure. Then we can define $\triangle^\times$ to be the closure of the set $\{\varphi \cdot \mathbb{P} : I_0 \subseteq I, \varphi = \prod_{i\in I_0} \varphi_i, \text{for } i \in I_0, \varphi_i \text{ is a conditional density from } \sigma\left(U_{J_i}\right) \text{ to } \sigma(U_i)\}$. Here a conditional density $\psi$ from $\sigma(X)$ to $\sigma(Y)$ is a $\sigma(X, Y)$-measurable probability density, s.t. $\mathbb{E}(\psi|X) = 1$. For example, a discrete Markov process fulfills this property with $J_i = \{i - 1\}$. It seems that the proofs in section Section 6 only us properties of $\triangle^\times$ that also hold in this case. Therefore it is highly likely that the whole theory generalizes to this case.

## 10.3 Continuous, ordered systems

In the previous section, we still had a discrete system. We could choose conditional densities independently. In continuous systems, we are not able to choose conditional densities independently. In a continuous time Markov process, the conditional probabilities are entangled. Let $\mathbb{P}$ be the law of a continuous time Markov process $M$ and $\varphi$ be a conditional density from time $M_t$ to $M_s$ Then $\varphi \cdot \mathbb{P}$ is no longer a Markov process. Therefore it is unclear how the theory can be generalized to this case.

## 10.4 Constraints

It could also be possible to extend the theory to allow for certain constraints on the probability distributions under considerations. For example, [12] introduces a set of probability distributions compatible with a hypergraph. This can be seen as putting a further constraint on $\triangle^{\times}$. More precisely, we let $\triangle^{\times}$ be the set of all product probability distributions $P$ that are absolutely continuous w.r.t. a reference measure $\mathbb{P}$, that also fulfill $N \in \mathcal{N} \Rightarrow P(N) = 0$ for some set system $\mathcal{N}$ that might be larger than the nullsets of $\mathbb{P}$. It is unclear if there is a generalization of the theory to this case.

## 10.5 Causal Discovery

We want to apply the theory to the problem of causal discovery. In Section 9, we have seen that structural independence corresponds to *d*-separation for causal structural models, and that it can be used in the toy example Example 9.3 to discover causal relations from observational data. It is likely that studying the relationship between structural independence and structural causal models will yield more examples of this kind. Ultimately, we would like to be able to axiomatize this relationship.

Furthermore, it seems likely that a theory for conditional systems with constraints could be used to infer more arrows of a graphical probabilistic system, if we allow the hypothesis class to include bayesnets with linear constraints.

## 10.6 Capturing all probability distributions

In this work, we have only considered product probability distributions that are mutually absolutely continuous, i.e. there is a reference measure $\mathbb{P}$. It is possible that the theory can be extended to all probability distributions, without a need for introducing a reference measure. Example 8.1 shows that the history cannot be a random index set in this case. To circumvent this problem, we could define an index set function as $J : \mathcal{A} \to \mathfrak{P}(I)$.

It seems possible that generation and disintegration generalize and give rise to a minimal generating index set function, which we would then call the history $\mathcal{H}(X|Z) : \mathcal{A} \to \mathfrak{P}(I)$. This definition would need to satisfy the fundamental theorem for all product probability distributions.

**Conjecture 10.6.1**: A history with type signature $\mathcal{A} \to \mathfrak{P}(I)$ exists and fulfills the fundamental theorem with all product probability distributions. Furthermore, the history w.r.t. a reference measure $\mathbb{P}$, introduced in this work, is a Radon-Nikodym derivative-like object w.r.t. $\mathbb{P}$ of the history as $\mathcal{A} \to \mathfrak{P}(I)$

# A Infinite product probability measures

This section is needed, to understand $\triangle^{\times}$ (c.f. Section 5) further. We first recall Kakutani's characterization of equivalent probability measures in countable product spaces [1], and apply it to our setting.

**In the following** Let $(\Omega, \mathcal{A}, \mathbb{P})$ be a probability space. Let $\triangle$ be the set of probability distributions on $(\Omega, \mathcal{A})$ that are absolutely continuous w.r.t. $\mathbb{P}$.

**Definition A.1**: There is an embedding $e : \triangle \to L^2(\Omega, \mathcal{A}, \mathbb{P})$, defined by $e(P) := \sqrt{\frac{\mathrm{d}P}{\mathrm{d}\mathbb{P}}}$. This embedding induces a Hilbert space structure on $\triangle$ with inner product $\langle P, Q\rangle = \int \sqrt{\frac{\mathrm{d}P}{\mathrm{d}\mathbb{P}} \frac{\mathrm{d}Q}{\mathrm{d}\mathbb{P}}} \,\mathrm{d}\mathbb{P} = \int \sqrt{\frac{\mathrm{d}P}{\mathrm{d}Q}} \,\mathrm{d}Q$. We denote the induced metric on $\triangle$ by $d_2 = d_2^{\mathbb{P}}$. For more details, see [1].

**Definition A.2**: We define the metric $d_1 = d_1^{\mathbb{P}}$ on $\triangle$ by $d_1(P, Q) := \int \left|\frac{\mathrm{d}P}{\mathrm{d}\mathbb{P}} - \frac{\mathrm{d}Q}{\mathrm{d}\mathbb{P}}\right| \mathrm{d}\mathbb{P}$.

**Lemma A.3**: $d_1(P, Q) = \int \left|\frac{\mathrm{d}P}{\mathrm{d}Q} - 1\right| \mathrm{d}Q$.

*Proof*: $d_1(P, Q) = \int \left|\frac{\mathrm{d}P}{\mathrm{d}\mathbb{P}} - \frac{\mathrm{d}Q}{\mathrm{d}\mathbb{P}}\right| \frac{\mathrm{d}\mathbb{P}}{\mathrm{d}Q} \,\mathrm{d}Q = \int \left|\frac{\mathrm{d}P}{\mathrm{d}\mathbb{P}} \frac{\mathrm{d}\mathbb{P}}{\mathrm{d}Q} - \frac{\mathrm{d}Q}{\mathrm{d}\mathbb{P}} \frac{\mathrm{d}\mathbb{P}}{\mathrm{d}Q}\right| \mathrm{d}Q = \int \left|\frac{\mathrm{d}P}{\mathrm{d}Q} - 1\right| \mathrm{d}Q$. □

**Definition A.4**: Measures $\mu$ and $\nu$ defined on the same measurable space are called equivalent, if they are mutually absolutely continuous.

**Theorem A.5**: Let $(\mu_n)_{n\in\mathbb{N}}$ and $(\nu_n)_{n\in\mathbb{N}}$ be families of probability measures. Then $\mu = \bigtimes_{n\in\mathbb{N}} \mu_n$ is equivalent to $\nu = \bigtimes_{n\in\mathbb{N}} \nu_n$ if and only if

$$\prod_{n\in\mathbb{N}} d_2(\mu_n, \nu_n) = \prod_{n\in\mathbb{N}} \int \sqrt{\frac{\mathrm{d}\mu_n}{\mathrm{d}\nu_n}} \,\mathrm{d}\nu_n > 0.$$

This condition is equivalent to the convergence of the series $\sum_{n\in\mathbb{N}} \log \int \sqrt{\frac{\mathrm{d}\mu_n}{\mathrm{d}\nu_n}} \,\mathrm{d}\nu_n$. In this case, we have $\prod_{n\in\mathbb{N}} \sqrt{\frac{\mathrm{d}\mu_n}{\mathrm{d}\nu_n}} \to \sqrt{\frac{\mathrm{d}\mu}{\mathrm{d}\nu}}$ in $L^2(\nu)$ and pointwise $\nu$-almost everywhere.

*Proof*: see [1]. □

**Lemma A.6**: $X_n \to X, Y_n \to Y$ in $L^2(\mathbb{P})$, Then $X_n Y_n \to XY$ in $L^1(\mathbb{P})$.

*Proof*: This is an immediate consequence of the Hölder inequality. □

**Corollary A.7**: If $P_n \overset{d_2}{\to} P$, then $P_n \overset{d_1}{\to} P$.

*Proof*: We have $\sqrt{\frac{\mathrm{d}P_n}{\mathrm{d}\mathbb{P}}} \overset{L^2(\mathbb{P})}{\to} \sqrt{\frac{\mathrm{d}P}{\mathrm{d}\mathbb{P}}}$. By Lemma A.6, $\frac{\mathrm{d}P_n}{\mathrm{d}\mathbb{P}} \overset{L^1(\mathbb{P})}{\to} \frac{\mathrm{d}P}{\mathrm{d}\mathbb{P}}$. □

The previous lemmas tell us that we can weaken the convergence in Theorem A.5 to $L^1$ or $d_1$ convergence.

We now apply these definitions and Kakutani's result to our setting with uncountably infinite products.

**Lemma A.8**: Let $P, Q \in \triangle^\times$ (c.f. Section 5). Then there is a family of positive densities $(\varphi_n)_{n\in\mathbb{N}_0}$ and a sequence of indices in $I$, $(i_n)_{n\in\mathbb{N}}$, s.t.

- $\mathbb{E}(\varphi_0 | U) \overset{\text{a.s.}}{=} 1$.
- $\forall n \in \mathbb{N} : \varphi_n$ is $\sigma\left(U_{i_n}\right)$-measurable.
- $\prod_{n\in\mathbb{N}_0} \varphi_n$ converges (unconditionally) in $L^1(P)$ and a.s. pointwise to $\frac{\mathrm{d}P}{\mathrm{d}Q}$.

*Proof*: Let $E$ denote the expectation w.r.t. $Q$. Let $\varphi = \frac{\mathrm{d}P}{\mathrm{d}Q}$, and set $\varphi_0 = \frac{\varphi}{E(\varphi|U)}$ and $\psi = E(\varphi|U)$. Then $\mathbb{E}(\varphi_0|U) = 1$, while $\varphi = \varphi_0 \psi$ and $\psi$ is $\sigma(U)$-measurable. By standard measure theory arguments, there exists $\psi' : \mathrm{Val}(U) \to \mathbb{R}$, s.t. $\psi'(U) = \psi$. We claim that $\psi' \cdot Q_U = P_U$. Let $A \subseteq \mathrm{Val}(U)$ measurable, then

$$\begin{aligned}
\psi' \cdot Q_U(A) &= \int_A \psi' \,\mathrm{d}Q_U = \int_{U^{-1}(A)} \psi'(U) \,\mathrm{d}Q = \int_{U^{-1}(A)} \psi \,\mathrm{d}Q = \int_{U^{-1}(A)} \psi \mathbb{E}(\varphi_0|U) \,\mathrm{d}Q \\
&= \int_{U^{-1}(A)} \mathbb{E}(\varphi|U) \,\mathrm{d}Q = \int_{U^{-1}(A)} \varphi \,\mathrm{d}Q = \int_{U^{-1}(A)} \mathrm{d}P = P_U(A).
\end{aligned}$$

Now $P_U$ and $Q_U$ are equivalent product probability measures on $\mathrm{Val}(U)$. More precisely, $P_U = \bigtimes_{i\in I} P_{U_i}$ and $Q_U = \bigtimes_{i\in I} Q_{U_i}$.

Let $\delta(i) := \log d_2\left(P_{U_i}, Q_{U_i}\right) = \log\left(\int \left(P_{U_i}/Q_{U_i}\right)^{\frac{1}{2}} \mathrm{d}Q_{U_i}\right)$. Clearly, for each $J \subseteq I$, we have $Q_{U_J}$ is equivalent to $P_{U_J}$. By Theorem A.5, we have $\sum_{i\in J} \delta_i < \infty$ for all $J \subseteq I$. Therefore, there can only be countably many $i$, s.t. $\delta_i \neq 0$. Define a sequence of indices s.t. $\{i_n : n \in \mathbb{N}\} \supseteq \{i \in I : \delta_i \neq 0\}$ and set $\varphi_n = \frac{\mathrm{d}P_{U_i}}{\mathrm{d}Q_{U_i}}(U_i)$. Since $P_U$ and $Q_U$ are equivalent product measures, we have $\prod_{n\in\mathbb{N}} \varphi_n = \frac{\mathrm{d}P_U}{\mathrm{d}Q_U}(U) = \psi$. Finally, $\varphi = \varphi_0 \psi = \prod_{n\in\mathbb{N}_0} \varphi_n$. The convergence properties of $\prod_{n\in\mathbb{N}} \varphi_n$ follow from Theorem A.5 or [1] and Lemma A.6. $\square$

Finally, it is important that a converging sequence of probability measures determines the limiting conditional expectation uniquely.

**Lemma A.9**: Let $P_n \overset{d_1}{\to} P$ and $X$ be bounded. Let $E$ and $E_n$ denote the expectation w.r.t $P$ and $P_n$. Then $E_n(X|Z) \to E(X|Z)$ in $P$-measure.

*Proof*: Since convergence in $P$-measure is metrizable, it suffices to show that any subsequence of $E_n(X|Z)$ has a subsequence that converges to $E(X|Z)$ in $P$-measure. W.l.o.g. it suffices to show that $E_n(X|Z)$ has a subsequence that converges to $E(X|Z)$ in $P$-measure. Let $\varphi_n$ be a positive density, s.t. $P_n = \varphi_n \cdot P$. Recall that $E_n(X|Z) = E(\varphi_n X|Z)/E(\varphi_n|Z)$. Since $\varphi_n \to 1$ in $L^1(P)$ (Lemma A.3) and $X$ is bounded, we have $E(\varphi_n|Z) \to \mathbb{E}(1|Z) = 1$ and $E(\varphi_n X|Z) \to E(X|Z)$ in $L^1(P)$. Therefore, by taking a non-relabeled subsequence we can assume that $E(\varphi_n|Z) \to 1$ and $E(\varphi_n X|Z) \to E(X|Z)$ $P$-a.s. Then clearly, $E_n(X|Z) \to E(X|Z)$ $P$-a.s. and thus in $P$-measure. $\square$

# B Well-known definitions and theorems

This section is intended to disambiguate notation. We don't provide proofs of the commonly known theorems. They can be found, for example, in [13]. Furthermore, we will apply the theorems without explicit reference to their appearance here. Let $\mathfrak{P}(\Omega)$ denote the powerset of $\Omega$. and for $S \subseteq \mathfrak{P}(\Omega)$, let $\sigma(S)$ denote the $\sigma$-algebra generated by $S$. Often, we assume that $\sigma(S)$ is completed w.r.t. a reference measure $\mathbb{P}$.

**Definition B.1** (measurable space): Let $\Omega$ be a set and $\mathcal{A} \subseteq \mathfrak{P}(\Omega)$ a $\sigma$-algebra. We call $(\Omega, \mathcal{A})$ a measurable space.

**Definition B.2** (random element): Let $(\Omega, \mathcal{A}), (\Omega', \mathcal{A}')$ be measurable spaces. We call a function $X : \Omega \to \Omega'$ that is $\mathcal{A}$-$\mathcal{A}'$-measurable a random element from $(\Omega, \mathcal{A})$ to $(\Omega', \mathcal{A}')$. It is convenient to introduce $X$ as a random element without referring explicitly to $(\Omega, \mathcal{A})$ or $(\Omega', \mathcal{A}')$, when $(\Omega, \mathcal{A})$ is understood from context. $(\Omega', \mathcal{A}')$ is denoted by $(\mathrm{Val}(X), \mathcal{V}(X))$.

**Definition B.3** (Dynkin system): A Dynkin system on $\Omega$ is a set $\mathcal{D} \subseteq \mathfrak{P}(\Omega)$, s.t. $\Omega \in \mathcal{D}$, and $A, B \in \mathcal{D} : A \subseteq B \Rightarrow B \setminus A \in \mathcal{D}$, and for $(A_n)_{n\in\mathbb{N}} \in \mathcal{D}^{\mathbb{N}}$ pairwise disjoint, $\bigcup_{n\in\mathbb{N}} A_n \in \mathcal{D}$.

**Definition B.4** ($\cap$-stable system): A $\cap$-stable system (or $\pi$-system), is a set $\mathcal{B} \subseteq \mathfrak{P}(\Omega)$, s.t. $A, B \in \mathcal{B} \Rightarrow A \cap B \in \mathcal{B}$.

**Theorem B.5**: If $S$ is a $\cap$-stable system and $\mathcal{D}$ a Dynkin system, then $S \subseteq \mathcal{D} \Rightarrow \sigma(S) \subseteq \mathcal{D}$.

**Definition B.6** (product measurable space): Given a family of measurable spaces $(\Omega_i, \mathcal{A}_i)_{i\in I}$, the product measurable space is defined by $\bigotimes_{i\in I} (\Omega_i, \mathcal{A}_i) := (\Omega, \mathcal{A})$, where $\Omega = \bigtimes_{i\in I}$ and $\mathcal{A}$ is the $\sigma$-algebra on $\Omega$ that renders the projections $\pi_i : \Omega \to \Omega_i$ measurable.

**Lemma B.7**: Given a product space $(\Omega, \mathcal{A}) = \bigotimes_{i \in I}(\Omega_i, \mathcal{A}_i)$, for all $A \in \mathcal{A}$, there is a countable set $I_0 \subseteq I$, s.t. $A \in \sigma(\pi_i : i \in I_0)$.

**Definition B.8**: Given a family of random elements $(X_i)_{i \in I}$ defined on a common measurable space $(\Omega, \mathcal{A})$, we associate to $(X_i)_{i \in I}$ the random element defined by $(\Omega, \mathcal{A}) \to \bigotimes_{i \in I}(\mathrm{Val}(X_i), \mathcal{V}(X_i)); \omega \mapsto (X_i(\omega))_{i \in I}$.

Let $(\Omega, \mathcal{A})$ be a measurable space and $\mathbb{P}$ be a probability distribution on this space. Let $\mathbb{E}$ denote the expectations of $\mathbb{P}$.

**Definition B.9** (density): A measurable function $\varphi : \Omega \to [0, \infty]$ is called density. It is called probability density w.r.t. $\mathbb{P}$, if $\mathbb{E}(\varphi) = 1$. We define the measure $\varphi \cdot \mathbb{P}$ by $(\varphi \cdot \mathbb{P})(A) \coloneqq \mathbb{E}(1_A \varphi)$.

**Definition B.10** (absolute continuity): We call $\mathbb{P}'$ absolutely continuous w.r.t. $\mathbb{P}$ and denote this by $\mathbb{P}' \ll \mathbb{P}$, if $\forall A \in \mathcal{A} : \mathbb{P}(A) = 0 \Rightarrow \mathbb{P}'(A) = 0$.

**Definition B.11** (mutual absolute continuity): We write $\mathbb{P} \sim \mathbb{P}'$, if $\mathbb{P} \ll \mathbb{P}'$ and $\mathbb{P}' \ll \mathbb{P}$.

**Theorem B.12** (Radon-Nikodym derivative): If $\mathbb{P}' \ll \mathbb{P}$, there exists a density $\varphi$ w.r.t. $\mathbb{P}$, s.t. $\varphi \cdot \mathbb{P} = \mathbb{P}'$. We write $\varphi =: \frac{\mathrm{d}\mathbb{P}'}{\mathrm{d}\mathbb{P}}$. If $\mathbb{P}' \sim \mathbb{P}$, then $\varphi$ can be chosen to be $\mathbb{P}$-a.s. positive.

**Definition B.13** (conditional expectation): Let $X \in L^1(\mathbb{P})$ and $\mathcal{C}$ a sub-$\sigma$-algebra of $\mathcal{A}$. There exists an up to $\mathbb{P}$-nullsets unique, $\mathcal{C}$-measurable map, $\mathbb{E}(X|\mathcal{C}) : \Omega \to \mathbb{R}$ that fulfills $\mathbb{E}(1_C X) = \mathbb{E}(1_C \mathbb{E}(X|Z))$ for all $C \in \mathcal{C}$. Let $Z$ be a random element. We set $\mathbb{E}(X|Z) = \mathbb{E}(X|\sigma(Z))$.

**Definition B.14** (conditional probability): Let $A \in \mathcal{A}$ and $Z$ be a random element. Then $\mathbb{P}(A|Z) \coloneqq \mathbb{P}(A|\sigma(Z)) \coloneqq \mathbb{E}(1_A|Z)$.

**Definition B.15** (conditional independence): Let $A, B \in \mathcal{A}$ and $\mathcal{C}$ a sub $\sigma$-algebra of $\mathcal{A}$. We say that $A$ is independent of $B$ given $\mathcal{C}$ w.r.t. $\mathbb{P}$, if $\mathbb{P}(A|\mathcal{C})\mathbb{P}(B|\mathcal{C}) \overset{\text{a.s.}}{=} \mathbb{P}(A, B|\mathcal{C})$. We write $A \perp\!\!\!\perp_{\mathbb{P}} B \mid \mathcal{C}$. We extend conditional independence to set systems. Let $\mathcal{A}_1, \mathcal{A}_2 \subseteq \mathcal{A}$. Then $\mathcal{A}_1 \perp\!\!\!\perp_{\mathbb{P}} \mathcal{A}_2 \mid \mathcal{C} :\Leftrightarrow \forall A \in \mathcal{A}_1, B \in \mathcal{A}_2 : A \perp\!\!\!\perp_{\mathbb{P}} B \mid \mathcal{C}$. We also allow the use of random elements. Let $X, Y, Z$ be random elements on $\Omega$. Then $X \perp\!\!\!\perp_{\mathbb{P}} Y \mid Z :\Leftrightarrow \sigma(X) \perp\!\!\!\perp_{\mathbb{P}} \sigma(Y) \mid \sigma(Z)$. Unconditional independence is written as $X \perp\!\!\!\perp_{\mathbb{P}} Y :\Leftrightarrow X \perp\!\!\!\perp_{\mathbb{P}} Y \mid \{\emptyset, \Omega\}$. Furthermore, independence given $Z$ holds for a family of random elements $(X_k)_{k \in K}$, if for all finite $K_0 \subseteq K$ and choices of $A_k \in \sigma(X_k)$ where $k \in K_0$, we have $\mathbb{P}\left(\bigcap_{k \in K} A_k | Z\right) = \prod_{k \in K} \mathbb{P}(A_k|Z)$. We use the word 'independence' to refer to both conditional and unconditional independence depending on the context.

**Lemma B.16**: Let $\mathcal{A}_1, \mathcal{A}_2 \subseteq \mathcal{A}$. If $\mathcal{A}_1 \perp\!\!\!\perp \mathcal{A}_2 \mid Z$, then this holds also for the smallest Dynkin systems $\mathcal{D}_1, \mathcal{D}_2$ larger than $\mathcal{A}_1, \mathcal{A}_2$, respectively. I.e. $\mathcal{D}_1 \perp\!\!\!\perp \mathcal{D}_2 \mid Z$.

*Proof sketch*: Let $A \in \mathcal{A}$. It is easy to see that $\{B \in \mathcal{B} : A \perp\!\!\!\perp B \mid Z\}$ forms a Dynkin system.□

In the following, let $X, Y$ and $Z$ be random elements on $(\Omega, \mathcal{A})$. Let $A, B \in \mathcal{A}$. Let $\mathbb{P}' \ll \mathbb{P}$ with expectation $\mathbb{E}'$. Let $\varphi$ be the density s.t. $\mathbb{P}' = \varphi \cdot \mathbb{P}$.

**Definition B.17**: $A \overset{\text{a.s.}}{=} B :\Leftrightarrow 1_A \overset{\text{a.s.}}{=} 1_B$.

**Lemma B.18**: $X \perp\!\!\!\perp_{\mathbb{P}} Y \mid Z \Leftrightarrow (X, Z) \perp\!\!\!\perp_{\mathbb{P}} (Y, Z) \mid Z$.

**Lemma B.19**: Let $A \in \mathcal{A}$. $A \perp\!\!\!\perp_{\mathbb{P}} Y \mid Z \Leftrightarrow \mathbb{P}(A|Y, Z) = \mathbb{P}(A|Z)$.

*Proof*: '$\Rightarrow$': Clearly, $\mathbb{P}(A|Z)$ is $\sigma(Y, Z)$-measurable. Let $C \in \sigma(Y, Z)$, then

$$\mathbb{E}(1_C \mathbb{P}(A|Z)) = \mathbb{E}(\mathbb{P}(A|Z)\mathbb{P}(C|Z)) = \mathbb{E}(\mathbb{P}(A \cap C|Z)) = \mathbb{P}(1_C 1_A).$$

'$\Leftarrow$': Let $B \in \sigma(Y)$ and $C \in \sigma(Z)$. Then

$\mathbb{E}(1_C \mathbb{P}(A|Z)\mathbb{P}(B|Z)) = \mathbb{E}(1_C \mathbb{E}(1_B \mathbb{P}(A|Z))) = \mathbb{E}(1_C 1_B \mathbb{P}(A|Y,Z)) = \mathbb{E}(1_C \mathbb{P}(A \cap B|Y,Z))$. $\square$

**Lemma B.20**: If $X \geq 0$ and $\mathbb{E}(X|Z) = 0$ then $X \overset{\text{a.s.}}{=} 0$.

*Proof*: $\mathbb{E}(X) = \mathbb{E}(\mathbb{E}(X|Z)) = 0$. $\square$

**Lemma B.21**: If $\mathbb{P}(A|Z) \overset{\text{a.s.}}{=} 1_A$, then $\exists C \in \sigma(Z)$ s.t. $C \overset{\text{a.s.}}{=} A$.

*Proof*: $C = \{\mathbb{P}(A|Z) > 0\} \in \sigma(Z)$. Then $C \overset{\text{a.s.}}{=} \{1_A > 0\}$. $\square$

**Lemma B.22**: $\mathbb{E}'(X|Z)\mathbb{E}(\varphi|Z) = \mathbb{E}(\varphi X|Z)$. If $\varphi > 0$, $\mathbb{E}'(X|Z) = \frac{\mathbb{E}(\varphi X|Z)}{\mathbb{E}(\varphi|Z)}$.

*Proof*: Clearly, the left hand side is $\sigma(Z)$-measurable. Let $C \in \sigma(Z)$. Then $\mathbb{E}(1_C \mathbb{E}'(X|Z)\mathbb{E}(\varphi|Z)) = \mathbb{E}(1_C \varphi \mathbb{E}'(X|Z)) = \mathbb{E}'(\mathbb{E}'(1_C X|Z)) = \mathbb{E}'(1_C X) = \mathbb{E}(1_C \varphi X)$. $\square$

**Lemma B.23**: If $\mathbb{P} \sim \mathbb{P}'$, then $A \perp\!\!\!\perp_{\mathbb{P}'} B \mid Z \Leftrightarrow \mathbb{E}(\varphi 1_A|Z)\mathbb{E}(\varphi 1_B|Z) = \mathbb{E}(\varphi|Z)\mathbb{E}(\varphi 1_A 1_B|Z)$ $\mathbb{P}$-a.s.

*Proof*: $\mathbb{P}'(A|Z)\mathbb{P}'(B|Z) = \mathbb{P}'(A,B|Z) \Leftrightarrow \mathbb{P}'(A|Z)\mathbb{E}(\varphi|Z)\mathbb{P}'(B|Z)\mathbb{E}(\varphi|Z) = \mathbb{P}'(A,B|Z)\mathbb{E}(\varphi|Z)$ $\Leftrightarrow \mathbb{P}(A|Z)\mathbb{P}(B|Z) = \mathbb{P}(A,B|Z)$. $\square$